# EXPONENTIAL GROWTH RATES IN A TYPED BRANCHING DIFFUSION

By Y. Git, J. W. Harris[1] and S. C. Harris

*Cambridge University, University of Bristol and University of Bath*

We study the high temperature phase of a family of typed branching diffusions initially studied in [*Astérisque* **236** (1996) 133–154] and [*Lecture Notes in Math.* **1729** (2000) 239–256 Springer, Berlin]. The primary aim is to establish some almost-sure limit results for the long-term behavior of this particle system, namely the speed at which the population of particles colonizes both space and type dimensions, as well as the rate at which the population grows within this asymptotic shape. Our approach will include identification of an explicit two-phase mechanism by which particles can build up in sufficient numbers with spatial positions near $-\gamma t$ and type positions near $\kappa\sqrt{t}$ at large times $t$. The proofs involve the application of a variety of martingale techniques—most importantly a "spine" construction involving a change of measure with an additive martingale. In addition to the model's intrinsic interest, the methodologies presented contain ideas that will adapt to other branching settings. We also briefly discuss applications to traveling wave solutions of an associated reaction–diffusion equation.

**1. Introduction.** In this article we will consider a certain family of typed branching diffusions that have particles which move (independently of each other) in space according to a Brownian motion with variance controlled by the particle's type process. The type of each particle evolves as an Ornstein–Uhlenbeck process and this type also controls the rate at which births occur. The particular form of this model permits many explicit calculations, but throughout we will strive to develop techniques that rely on general principles as much as possible, so they might readily adapt to other situations. This model was previously considered in [12, 13]; these papers form essential foundations for this work, although we will recall various results as necessary.

Received December 2004; revised November 2006.
[1]Supported in part by an EPSRC studentship.

We will make some significant applications of the spine theory for branching processes. Inspired by the series of papers Lyons, Pemantle and Peres [19], Lyons [18] and Kurtz, Lyons, Pemantle and Peres [16], spine techniques have been instrumental in recent years in providing intuitive and elegant proofs of many important classical and new results in the theory of branching processes. In this article we use the recent reformulation of the spine method presented in [8], which follows in similar spirit to the branching Brownian motion study of Kyprianou [17]. For a selection of other applications of spine techniques, for example, see [1, 6, 7, 23] and references therein.

1.1. *The branching model.* We define $N_t$ to be the set of particles alive at time $t \geq 0$. For a particle $u \in N_t$, $X_u(t) \in \mathbb{R}$ is its spatial position, and $Y_u(t) \in \mathbb{R}$ is the *type* of $u$. We will label offspring using the Ulam–Harris convention where, for example, if $u = \varnothing 21$ then particle $u$ is the first child of the second child of the initial ancestor, and we will write $v > u$ if particle $v$ is a descendant of particle $u$. The configuration of the branching diffusion at time $t$ is given by the point process $\mathbb{X}_t := \{(X_u(t), Y_u(t)) : u \in N_t\}$.

A particle's type evolves as an Ornstein–Uhlenbeck process with an invariant measure given by the standard normal density $\phi(y)$ and an associated differential operator (generator)

$$\mathcal{Q}_\theta := \frac{\theta}{2}\left(\frac{\partial^2}{\partial y^2} - y\frac{\partial}{\partial y}\right),$$

where $\theta > 0$ is considered to be the *temperature* of the system. The spatial motion of a particle of type $y$ is a driftless Brownian motion on $\mathbb{R}$ with variance

$$A(y) := ay^2, \qquad \text{where } a \geq 0.$$

A particle of type $y$ particle is replaced by two offspring at a rate

$$R(y) := ry^2 + \rho, \qquad \text{where } r, \rho \geq 0.$$

Each offspring inherits its parent's current type and spatial position, and then moves off independently of all others. We use $P^{x,y}$ and $E^{x,y}$ with $x, y \in \mathbb{R}$ to represent probability and expectation when the Markov process starts with a single particle at position $(x, y)$.

We will find the almost sure rate of exponential growth, $D(\gamma, \kappa)$, of particles which are found simultaneously with spatial positions near $-\gamma t$ and type positions near $\kappa\sqrt{t}$ at large times $t$. From this we can deduce the speed of extremal particles and hence the asymptotic shape of the particle system. The main effort is required in identifying $D(\gamma, \kappa)$ as the almost sure limit of

$$t^{-1} \log \sum_{u \in N_t} \mathbb{1}_{\{X_u(t) \leq -\gamma t; Y_u(t) \geq \kappa\sqrt{t}\}}.$$

In particular, the convergence properties of two different families of additive martingales associated with the branching diffusion will lead directly to the spatial exponential growth rates and an upper bound on the space-type growth. For the remaining lower bound, we describe an explicit two-phase mechanism for amassing the required number of particles with prescribed space-type positions. The first phase involves building up an "excess" number of particles, each covering a certain proportion of the required spatial distance. During their second phase, enough of these particles must succeed in making a difficult and rapid ascent into the required position. The latter phase is proved using an intuitive change of measure technique that induces a spine construction.

The family of models we are considering is specific but nevertheless have some features of fundamental significance that motivate the choices for $\mathcal{Q}_\theta$, $R$ and $A$. If the spatial motion is ignored, we have investigated a binary branching Ornstein–Uhlenbeck process in a quadratic breeding potential. In contrast, Enderle and Hering [5] considered a branching Ornstein–Uhlenbeck with constant branching rate but random offspring distribution. A quadratic breeding potential is a critical rate for *explosions* in the population of particles. In a branching Brownian motion on $\mathbb{R}$ with binary splitting occurring at rate $x^p$ at position $x$, the population will explode almost surely in finite time if $p > 2$, whereas for $p = 2$ the expected number of particles explodes while the total population remains finite for all time with probability 1 (see [15], Chapter 5.12). The Ornstein–Uhlenbeck process is not only a canonical ergodic diffusion, but this type-motion has exactly the right drift to help counteract the quadratic breeding rate. For high temperatures, $\theta > 8r$, there is a sufficiently strong mean-reversion in the type processes to ensure that the expected total population size does not blow up; but for temperatures $\theta \leq 8r$, the quadratic breeding overpowers the pull toward the origin, the expected population blows up in a finite time and particles behave very differently. Throughout this paper we consider only high temperatures $\theta > 8r$, deferring the low and critical temperature regimes to future work. Given other choices, the quadratic spatial diffusion coefficient now becomes very natural, enabling us to find *explicit* families of (fundamental) additive martingales since the linearized traveling-wave equation can be linked to the classical harmonic oscillator equations from physics. The binary branching mechanism was taken for simplicity; in principle our approach could extend to general offspring distributions, although new features would arise from possible extinctions and necessary offspring moment conditions. All these choices make the models rich in structure, possessing some very challenging features whilst remaining sufficiently tractable.

1.2. *Application to reaction–diffusion equations.* Following in the footsteps of McKean [20], the solution of the reaction–diffusion equation

$$\frac{\partial u}{\partial t} = \frac{1}{2}A(y)\frac{\partial^2 u}{\partial x^2} + R(y)u(u-1) + \frac{\theta}{2}\left(\frac{\partial^2 u}{\partial y^2} - y\frac{\partial u}{\partial y}\right), \quad (1)$$

with initial condition $f(x,y) \in [0,1]$ for all $x,y \in \mathbb{R}$, can be represented by

$$u(t,x,y) = E^{x,y}\left(\prod_{u \in N_t} f(X_u(t), Y_u(t))\right). \quad (2)$$

Of great importance for reaction–diffusion equations are traveling-wave solutions (e.g., see [21]). In the present context, a solution to equation (1) of the form $u(t,x,y) := w(x-ct,y)$ is said to be a *traveling-wave of speed c*, where $w(x,y)$ solves the traveling-wave equation

$$\frac{1}{2}A(y)\frac{\partial^2 w}{\partial x^2} + c\frac{\partial w}{\partial x} + R(y)w(w-1) + \frac{\theta}{2}\left(\frac{\partial^2 w}{\partial y^2} - y\frac{\partial w}{\partial y}\right) = 0. \quad (3)$$

Fundamental to our study of the branching diffusion are two families of "additive" martingales, $Z_\lambda^\pm(t)$ [defined at (6)], which are linked to the linearization of (1). When $\theta > 8r$, Harris and Williams [13] determined when $Z_\lambda^-$ is uniformly integrable (see Theorem 17) and then $w_\lambda(x,y) := E^{x,y}\exp(-Z_\lambda^-(\infty))$ yields a traveling wave of speed $c_\lambda^-$. This gives the existence of traveling waves for all speeds $c$ greater than some threshold $\tilde{c}(\theta) := \inf c_\lambda^-$.

Furthermore, combining the McKean representation (2) with the almost-sure convergence result established in [12] (look ahead to Theorem 18) can give results on the attraction toward traveling waves from given initial data. For example, if $-\ln f(x,y) \sim e^{\lambda x}g(y)$ *uniformly* in $y$ as $x \to \infty$ for some suitable $g \in L^2(\phi)$, the solution $u(x,y)$ to (1) with initial conditions $f$ satisfies $u(t, x - c_\lambda^- t, y) \to w_\lambda(x + \hat{x}, y)$ as $t \to \infty$, where $\hat{x}$ is some constant that can be determined from $g$.

In future work we hope to develop the approach used for standard BBM and the FKPP equation in [11], and prove that traveling waves of a given speed $c > \tilde{c}(\theta)$ are *unique* (up to translation) and that no traveling waves exist for speeds $c < c(\theta)$. We anticipate that our new results on the growth rates of particles will aid in establishing some difficult estimates on the tail behavior of any traveling wave, and hence assist in proving the conjectured uniqueness. In addition, we expect our growth rate results will be essential in obtaining broader classes of initial conditions that are attracted toward traveling waves. In each of these problems, difficulties arise from the unbounded type space where, for example, some control must be gained over the possible contributions to $\sum_{u \in N_t} \log f(X_u(t) - ct, Y_u(t))$ from particles that have large type positions in addition to large spatial positions.

**2. Main results.** In this section, we will present our main results that identify the growth rates found within the branching diffusion. We will give an overview for our proofs, identifying the key ideas and techniques used, as well as introducing some intuition for the dominant behavior of particles that underpins our approach.

2.1. *Martingales.* The principal tools used throughout this paper are two fundamental families of "additive" martingales, which were introduced in [13].

Before defining the martingales we give some key definitions. Let

$$\lambda_{\min} := -\sqrt{\frac{\theta - 8r}{4a}}.$$

Let $\lambda \in \mathbb{R}$, with the following convention which we always use for $\lambda$:

$$\lambda_{\min} < \lambda < 0.$$

Also, define

(4) $$\mu_\lambda := \frac{1}{2}\sqrt{\theta(\theta - 8r - 4a\lambda^2)}, \qquad \psi_\lambda^\pm := \frac{1}{4} \pm \frac{\mu_\lambda}{2\theta},$$

(5) $$E_\lambda^\pm := \rho + \theta \psi_\lambda^\pm, \qquad c_\lambda^\pm := -E_\lambda^\pm/\lambda.$$

Will will occasionally write $E_\lambda^\pm$ as $E^\pm(\lambda)$ in order to emphasize that $E_\lambda^\pm$ are really functions of $\lambda$; the $\pm$ superscripts will always distinguish these from expectation operators. Note that $\lambda_{\min}$ is the point beyond which $\mu_\lambda$ is no longer a real number.

The martingales are $Z_\lambda^-$ and $Z_\lambda^+$, defined for $\lambda \in (\lambda_{\min}, 0]$ as

(6) $$Z_\lambda^\pm(t) := \sum_{u \in N_t} v_\lambda^\pm(Y_u(t))e^{\lambda X_u(t) - E_\lambda^\pm t},$$

where $v_\lambda^\pm(y) := \exp(\psi_\lambda^\pm y^2)$ are strictly-positive eigenfunctions of the operator

$$\mathcal{Q}_\theta + \tfrac{1}{2}\lambda^2 A + R,$$

with corresponding eigenvalues $E_\lambda^- < E_\lambda^+$ and $A, R$ are the functions defined in Section 1.1. This operator is self-adjoint on $L^2(\phi)$ with the inner product $\langle \cdot, \cdot \rangle_\phi$ where $\langle f, g \rangle_\phi := \int fg\phi \, dy$ and $\phi$ is the standard normal density. Note that $v_\lambda^- \in L^2(\phi)$, whereas $v_\lambda^+ \notin L^2(\phi)$ so is *not* normalizable.

The calculations of Section 3 make it easy to see these are martingales, and throughout the paper we will need a variety of martingale convergence results which are gathered together in Section 8. In particular, we will need to know precisely when $Z_\lambda^-$ is uniformly integrable with a strictly positive limit, some further strong convergence results for other closely related sums over particles (also identifying which particles contribute nontrivially to their limits), and the rate of convergence to zero of the $Z_\lambda^+$ martingales.

2.2. *The asymptotic growth-rate of particles along spatial rays.* As an essential initial step toward determining the growth rate of particles in the two-dimensional space-type domain, we first look at the growth rate of particles in the spatial dimension only.

For $\gamma \geq 0$ and $C \subset \mathbb{R}$, define

$$N_t(\gamma; C) := \sum_{u \in N_t} \mathbb{1}_{\{X_u(t) \leq -\gamma t; Y_u(t) \in C\}}. \tag{7}$$

The limit giving the *expected* rate of growth,

$$\lim_{t \to \infty} t^{-1} \log E(N_t(\gamma; \mathbb{R}))$$

can be shown to exist and its value can be calculated to be

$$\begin{aligned}\Delta(\gamma) &:= \inf_{\lambda \in (\lambda_{\min}, 0)} \{E_\lambda^- + \lambda \gamma\} \\ &= \rho + \frac{\theta}{4} - \frac{1}{4}\sqrt{a^{-1}(\theta - 8r)(4\gamma^2 + \theta a)}.\end{aligned} \tag{8}$$

An outline for this expectation calculation is given in Section 3.

It is now tempting to guess that the asymptotic speed of the spatially left-most particle, $\tilde{c}(\theta)$, is given by

$$\begin{aligned}\tilde{c}(\theta) &:= \sup\{\gamma : \Delta(\gamma) > 0\} \\ &= \sqrt{2a\left(r + \rho + \frac{2(2r + \rho)^2}{\theta - 8r}\right)}.\end{aligned} \tag{9}$$

Recall that $\tilde{c}(\theta) = \inf_{\lambda \in (\lambda_{\min}, 0)} c_\lambda^-$ is also the minimum threshold for traveling waves. In this particular situation, the guess that "expectation" and "almost sure" right-most particle speeds agree was first proved rigorously using a martingale change of measure technique in [13]. In this paper, we extend this connection and prove that the "expected" and "almost sure" *rates of growth* of particles with given speeds (Theorem 1) and given space-type locations (Theorem 3) agree.

THEOREM 1. *Let $\gamma \geq 0$ and $y_0 < y_1$. Under each $P^{x,y}$ law, the limit*

$$D(\gamma) := \lim_{t \to \infty} t^{-1} \log N_t(\gamma; [y_0, y_1])$$

*exists almost surely and is given by*

$$D(\gamma) = \begin{cases} \Delta(\gamma), & \text{if } 0 \leq \gamma < \tilde{c}(\theta), \\ -\infty, & \text{if } \gamma \geq \tilde{c}(\theta). \end{cases}$$

Note that symmetry in the process means there is a corresponding result for particles with spatial velocities greater than $+\gamma$ (corresponding to positive $\lambda$ values). We may occasionally make use of such process symmetries without further comment. Then, since $N_t(\gamma;\mathbb{R})$ is integer valued, the asymptotic speed of the right-most particle follows immediately:

COROLLARY 2. *Almost surely,*
$$\lim_{t\to\infty} t^{-1} \sup\{X_u(t) : u \in N_t\} = \tilde{c}(\theta).$$

This spatial growth rate result is proved in Section 10 using the martingale results from Section 8. In fact, it is very easy to obtain the upper bound by first dominating the indicator function with exponentials to reveal that $N_t(\gamma;\mathbb{R}) \leq \exp\{(E_\lambda^- + \lambda\gamma)t\}Z_\lambda^-(t)$, recalling that $Z_\lambda^-$ is a convergent martingale, and then optimizing over the choice of $\lambda$. For the lower bound, we will use a strong convergence result obtained in [12], combined with the idea that each uniformly integrable martingale $Z_\lambda^-$ essentially "counts" only the particles of corresponding velocity $-\gamma$.

2.3. *The asymptotic shape and growth of the branching diffusion.* The main result of this paper is the almost-sure rate of growth of particles which are in the vicinity of $-\gamma t$ in space and near $\kappa\sqrt{t}$ in type position at large times $t$. For $\gamma, \kappa \geq 0$, it can be shown that the limit

(10) $$\lim_{t\to\infty} t^{-1} \log E(N_t(\gamma; [\kappa\sqrt{t}, \infty)))$$

exists and takes the value

(11) $$\Delta(\gamma,\kappa) := \inf_{\lambda \in (\lambda_{\min},0)} \{E_\lambda^- + \lambda\gamma - \kappa^2\psi_\lambda^+\}$$
$$= \rho + \frac{(\theta - \kappa^2)}{4} - \frac{1}{4\theta a}\sqrt{\theta(\theta-8r)(4a\theta\gamma^2 + a^2(\theta+\kappa^2)^2)}.$$

An outline of this expectation calculation is given in Section 3. Once again, we will find that the "almost sure" rate of growth of particles agrees with this "expected" rate exactly where there is *growth* in particle numbers.

THEOREM 3. *Let $\gamma, \kappa \geq 0$ with $\Delta(\gamma,\kappa) \neq 0$. Under each $P^{x,y}$ law, the limit*
$$D(\gamma,\kappa) := \lim_{t\to\infty} t^{-1} \log N_t(\gamma; [\kappa\sqrt{t}, \infty))$$

*exists almost surely and is given by*

(12) $$D(\gamma,\kappa) = \begin{cases} \Delta(\gamma,\kappa), & \text{if } \Delta(\gamma,\kappa) \geq 0, \\ -\infty, & \text{if } \Delta(\gamma,\kappa) < 0. \end{cases}$$

To prove the tricky lower bound of Theorem 3, which amounts to the major work of this paper, we will exhibit an explicit two-phase mechanism by which the branching diffusion can build up *at least* the required exponential number of particles near to $-\gamma t$ in space and $\kappa\sqrt{t}$ in type position by large times $t$.

During the first phase, over a large time $t$ the process builds up an initial excess of approximately $\exp(\Delta(\alpha)t)$ particles with spatial position at least $-\alpha t$, as is already known from Theorem 1. In this "ergodic" phase, "typical" particles found near $-\alpha t$ in space will have drifted with a steady spatial speed of $\alpha$ whilst their type histories will have behaved roughly like OU processes with *inward* drift of $\mu_\lambda y$ for a certain optimal choice $\lambda(\alpha)$ of parameter $\lambda$.

For the second phase, we will show that the probability any individual particle has *at least one descendant* that makes a "rapid ascent" in both space and type dimensions from initial position $(0,0)$ to final position near $(-\beta t, \kappa\sqrt{t})$ is approximately $\exp(-\Theta(\beta,\kappa)t)$, where

$$\text{(13)} \qquad \Theta(\beta,\kappa) = \frac{\kappa^2}{4} + \frac{\sqrt{\theta(\theta-8r)(a^2\kappa^4 + 4a\theta\beta^2)}}{4a\theta},$$

and the time taken for this "rapid ascent" is an interval $[0,\tau]$. We show that this time $\tau$ can be chosen such that $2\mu_\lambda\tau \sim \log t$, and hence the additional time is asymptotically negligible in comparison with $t$. Intuitively, we will see that *given* an offspring that has successfully made such a difficult "rapid ascent," it will most likely have roughly had its type process behaving like an OU process with an *outward* drift of $\mu_\lambda y$ and the Brownian motion driving its spatial motion will have had a drift $\lambda$ [corresponding to a real time spatial drift $\lambda A(y)$ that increases in strength as the type position $y$ increases], for some optimal choice $\lambda(\beta,\kappa)$ of parameter $\lambda$. The precise result required will be formulated rigorously as a large-deviation lower bound in Theorem 7 of Section 5, and is proved using a "spine" change of measure technique intimately related to the $Z_\lambda^+$ martingales.

Combining these two phases and using independence of the particles, we can see that the number of particles near $(-\alpha t, 0)$ at time $t$ that subsequently proceed to have *at least one* descendant near $(-(\alpha+\beta)t, \kappa\sqrt{t})$ is *approximately* Poisson with mean

$$\exp(\{\Delta(\alpha) - \Theta(\beta,\kappa)\}t).$$

Optimizing for a fixed overall spatial speed $\gamma$, some calculus reveals that

$$\text{(14)} \qquad \sup_{\substack{\alpha+\beta=\gamma \\ \alpha,\beta>0}} \{\Delta(\alpha) - \Theta(\beta,\kappa)\} = \Delta(\bar{\alpha}) - \Theta(\bar{\beta},\kappa) = \Delta(\gamma,\kappa),$$

with optimal parameters

(15) $$\bar{\alpha} = \gamma\left(\frac{\theta}{\theta + \kappa^2}\right) \quad \text{and} \quad \bar{\beta} = \gamma\left(\frac{\kappa^2}{\theta + \kappa^2}\right).$$

Thus we will be able to demonstrate an explicit two-phase mechanism producing the required number of particles, with this outline argument later guiding our rigorous proof. In addition, it is interesting to note that the optimal choices for $\lambda$ over each phase then also coincide at a single value $\bar{\lambda} = \lambda(\bar{\alpha}) = \lambda(\bar{\beta}, \kappa)$.

An informative large deviation *heuristic* for the rapid ascent can also be found in Section 4, with this section also containing some essential optimal path calculations. We actually prove the two-phase mechanism for the lower bound of Theorem 3 in Section 5, although we defer proving the large-deviation lower bound until Section 7 after presenting the necessary "spine" background in Section 6.

We prove the upper bound for the space-type growth rate in Section 9, again making crucial use of martingale results from Section 8. Similarly to the spatial growth case, we can find an upper bound using the $Z_\lambda^+$ martingales, that is to say $N_t(\gamma; [\kappa\sqrt{t}, \infty)) \leq \exp\{(E_\lambda^+ + \lambda\gamma - \kappa^2\psi_\lambda^+)t\}Z_\lambda^+(t)$. However, as each $Z_\lambda^+$ martingale converges to *zero*, we must show that its exponential decay rate is $(E_\lambda^+ - E_\lambda^-)$ before being able to optimize over the choice of $\lambda$ to obtain the required upper bound.

Given Theorem 3, and noting symmetries, it becomes straightforward to retrieve the following:

COROLLARY 4. *For any $F \subset \mathbb{R}^2$, define*

$$\mathcal{N}_t(F) := \sum_{u \in N_t} \mathbb{1}_{\{(X_u(t)/t, Y_u(t)/\sqrt{t}) \in F\}}.$$

*If $B \subset \mathbb{R}^2$ is any open set and $C \subset \mathbb{R}^2$ is any closed set, then almost surely under any $P^{x,y}$*

$$\liminf_{t \to \infty} \frac{1}{t} \log \mathcal{N}_t(B) \geq \sup_{(\gamma,\kappa) \in B} D(\gamma, \kappa),$$

$$\limsup_{t \to \infty} \frac{1}{t} \log \mathcal{N}_t(C) \leq \sup_{(\gamma,\kappa) \in C} D(\gamma, \kappa),$$

*with the growth rate $D(\gamma, \kappa)$ given at equation* (12).

We can also recover the almost sure asymptotic shape of the region occupied by the particles in the branching diffusion.

COROLLARY 5. *Let $B \subset \mathbb{R}^2$ be any open set. Almost surely, under each $P^{x,y}$ law,*

$$\mathcal{N}_t(B) \to \begin{cases} 0, & \text{if } \mathcal{S} \cap B = \varnothing, \\ +\infty, & \text{if } \mathcal{S} \cap B \neq \varnothing, \end{cases}$$

*where $\mathcal{S} \subset \mathbb{R}^2$ is the set given by*

$$\mathcal{S} := \{(\gamma, \kappa) \in \mathbb{R}^2 | \Delta(\gamma, \kappa) > 0\}.$$

**3. Some expectation calculations.** This section discusses how the *expected* growth rates given in the previous section may be obtained. For this, we use the "many-to-one" lemma (see, e.g., [8]) and one-particle changes of measure. In the process we shall start to gain valuable intuition into how particles within the branching diffusion behave, as well as seeing hints as to which are the "correct" martingales to use to prove the almost-sure growth rate results.

For simplicity, we assume throughout this section that the branching diffusion starts with one particle at the origin in both space and type at time zero, unless otherwise stated. We also introduce a family of single particle probability measures $\mathbb{P}_{\mu,\lambda}$ with associated expectations $\mathbb{E}_{\mu,\lambda}$ where, under $\mathbb{P}_{\mu,\lambda}$, $\eta$ is an Ornstein–Uhlenbeck process with variance $\theta$ and drift $\mu$, and $\xi_t = B(\int_0^t A(\eta_s)\,ds)$ where $B$ is a Brownian motion with drift $\lambda$.

LEMMA 6 (Many-to-one). *If $f : \mathbb{R}^2 \mapsto \mathbb{R}$ is Borel measurable then*

$$(16) \quad E \sum_{u \in N_t} f(X_u(t), Y_u(t)) = \mathbb{E}_{\theta/2, 0}\left(\exp\left(\int_0^t R(\eta_s)\,ds\right) f(\xi_t, \eta_t)\right).$$

Using the many-to-one lemma, and changing measure to alter the drift of Brownian motion, we see that

$$E \sum_{u \in N_t} f(X_u(t), Y_u(t))$$

$$= \mathbb{E}_{\theta/2, 0}\left(\exp\left(\int_0^t R(\eta_s)\,ds\right) f(\xi_t, \eta_t)\right)$$

$$= \mathbb{E}_{\theta/2, 0}\left(e^{-\lambda \xi_t} \exp\left(\int_0^t \left\{R(\eta_s) + \frac{\lambda^2}{2} A(\eta_s)\right\} ds\right)\right.$$

$$\left. \times f(\xi_t, \eta_t) \cdot e^{\lambda \xi_t - \lambda^2/2 \int_0^t A(\eta_s)\,ds}\right)$$

$$= \mathbb{E}_{\theta/2, \lambda}\left(\exp\left(-\lambda \xi_t + \int_0^t \left\{R(\eta_s) + \frac{\lambda^2}{2} A(\eta_s)\right\} ds\right) f(\xi_t, \eta_t)\right).$$

To perform a further change of measure on the OU process to get rid of the time integrals in the exponential of the expectation, we recall that

$$\left.\frac{d\mathbb{P}_{\mu_\lambda,\cdot}}{d\mathbb{P}_{\theta/2,\cdot}}\right|_{\mathcal{F}_t} = M_t^{\mu_\lambda,\theta/2}$$

$$:= \exp\left(\psi_\lambda^- \eta_t^2 - E_\lambda^- t + \int_0^t \left\{R(\eta_s) + \frac{1}{2}\lambda^2 A(\eta_s)\right\} ds\right)$$

and then

$$E \sum_{u \in N_t} f(X_u(t), Y_u(t))$$

(17)
$$= \mathbb{E}_{\theta/2,\lambda}(\exp(-\lambda \xi_t - \psi_\lambda^- \eta_t^2 + E_\lambda^- t) f(\xi_t, \eta_t) \cdot M_t^{\mu_\lambda,\theta/2})$$

$$= \mathbb{E}_{\mu_\lambda,\lambda}(\exp(-\lambda \xi_t - \psi_\lambda^- \eta_t^2 + E_\lambda^- t) f(\xi_t, \eta_t)).$$

Note that the many-to-one lemma, combined with the branching property, immediately suggests how to get "additive" martingales for the branching diffusion from single particle martingales—for example, taking $f(x,y) = \exp\{\lambda x + \psi_\lambda^- y^2\}$ in equation (17) quickly leads to the martingale $Z_\lambda^-$.

We may now proceed to calculate the *expected* growth rates. However, for both clarity and brevity we will leave rigorous details to the interested reader, noting that the intuition we will gain from our rough calculations will later be invaluable in guiding our rigorous proof of the corresponding *almost-sure* growth rates.

3.1. *The expected rate of growth along spatial rays.* We first give the outline of some calculations to find the rate of growth in the expected number of particles near $-\gamma t$ in space at time $t$.

Using the formula from (17), for $\lambda \in (\lambda_{\min}, 0)$ and any $\varepsilon > 0$ we have

$$E \sum_{u \in N_t} \mathbb{1}_{\{t^{-1} X_u(t) + \gamma \in (-\varepsilon, \varepsilon)\}}$$

$$= \mathbb{E}_{\mu_\lambda,\lambda}(e^{-\lambda \xi_t - \psi_\lambda^- \eta_t^2 + E_\lambda^- t} \mathbb{1}_{\{t^{-1} \xi_t + \gamma \in (-\varepsilon, \varepsilon)\}})$$

$$\begin{cases} \leq e^{(E_\lambda^- + \lambda\gamma - \lambda\varepsilon)t} \mathbb{E}_{\mu_\lambda,\lambda}\left(e^{-\psi_\lambda^- \eta_t^2}; \frac{\xi_t}{t} + \gamma \in (-\varepsilon, \varepsilon)\right) \\ \geq e^{(E_\lambda^- + \lambda\gamma + \lambda\varepsilon)t} \mathbb{E}_{\mu_\lambda,\lambda}\left(e^{-\psi_\lambda^- \eta_t^2}; \frac{\xi_t}{t} + \gamma \in (-\varepsilon, \varepsilon)\right) \end{cases}$$

where, with some abuse of notation that we shall continue to use throughout this section, we will abbreviate this to

(18)
$$E \sum_{u \in N_t} \mathbb{1}_{\{X_u(t) \sim -\gamma t\}} = \mathbb{E}_{\mu_\lambda,\lambda}(e^{-\lambda \xi_t - \psi_\lambda^- \eta_t^2 + E_\lambda^- t} \mathbb{1}_{\{\xi_t \sim -\gamma t\}})$$

$$\sim e^{(E_\lambda^- + \lambda\gamma)t} \mathbb{E}_{\mu_\lambda,\lambda}(e^{-\psi_\lambda^- \eta_t^2}; \xi_t \sim -\gamma t)$$

with the understanding that any subsequent arguments to identify exponential growth rates can readily be made rigorous by using the appropriate upper and lower bounds, and so on.

Now, considering $E^-(\lambda) := E^-_\lambda$ as a function of $\lambda$, we have from (8) that $\Delta(\gamma) = \inf_{\lambda \in (\lambda_{\min}, 0)} \{E^-(\lambda) + \lambda\gamma\} = E^-(\lambda_\gamma) + \lambda_\gamma \gamma$, where $\lambda_\gamma$ satisfies

$$(19) \qquad \frac{\partial E^-}{\partial \lambda}(\lambda_\gamma) = -\gamma, \qquad \text{hence } \lambda_\gamma = -\gamma\sqrt{\frac{(\theta - 8r)}{\theta a^2 + 4a\gamma^2}}.$$

Of course, choosing this optimal $\lambda_\gamma$ value in (18) means that we must have simultaneously maximized the expectation $\mathbb{E}_{\mu_\lambda, \lambda}(\exp(-\psi^-_\lambda \eta_t^2); \xi_t \sim -\gamma t)$, and to confirm that this value is not exponentially decaying in $t$ is now relatively straightforward. Under $\mathbb{P}_{\mu_\lambda, \lambda}$, $\eta$ is an Ornstein–Uhlenbeck process with an invariant measure given by the probability density, $\phi_\lambda$, of the normal distribution $N(0, \theta/(2\mu_\lambda))$; and $\xi_t = B(\int_0^t A(\eta_s) \, ds)$, where $B$ is a BM with drift $\lambda$. Note also that by differentiating

$$\langle (\mathcal{Q}_\theta + (1/2)\lambda^2 A + R - E^-_\lambda)v^-_\lambda, v^-_\lambda \rangle_\phi = 0$$

with respect to $\lambda$, using self-adjointness, and observing that $\phi_\lambda \propto (v^-_\lambda)^2 \phi$, we find that

$$\frac{\partial E^-}{\partial \lambda} = \frac{\langle \lambda A v^-_\lambda, v^-_\lambda \rangle_\phi}{\langle v^-_\lambda, v^-_\lambda \rangle_\phi} = \lambda \int_\mathbb{R} A(y) \phi_\lambda(y) \, dy.$$

Then almost surely under $\mathbb{P}_{\mu_\lambda, \lambda}$,

$$(20) \quad \frac{\xi_t}{t} = \frac{B(\int_0^t A(\eta_s) \, ds)}{\int_0^t A(\eta_s) \, ds} \frac{\int_0^t A(\eta_s) \, ds}{t} \to \lambda \int_\mathbb{R} A \phi_\lambda \, dy = \frac{\lambda a \theta}{2\mu_\lambda} = \frac{\partial E^-}{\partial \lambda},$$

and so when we use the optimal $\lambda_\gamma$ value we get exactly the desired drift, since $\frac{\partial E^-}{\partial \lambda}(\lambda_\gamma) = -\gamma$. Then

$$\mathbb{E}_{\mu_{\lambda_\gamma}, \lambda_\gamma}(e^{-\psi^-_{\lambda_\gamma} \eta_t^2}; \xi_t \sim -\gamma t) \to \lim_{t \to \infty} \mathbb{E}_{\mu_{\lambda_\gamma}, \lambda_\gamma}(e^{-\psi^-_{\lambda_\gamma} \eta_t^2})$$

$$= \int_\mathbb{R} e^{-\psi^-_{\lambda_\gamma} y^2} \phi_{\lambda_\gamma}(y) \, dy.$$

In this way, we can obtain the exact rate of exponential growth for the expectation,

$$\lim_{t \to \infty} t^{-1} \log E(N_t(\gamma; \mathbb{R})) = \Delta(\gamma).$$

The changes of measure used above are actually suggesting a great deal about the dominant particles that are found in the vicinity of a given ray in space. An alternative discussion of this expectation result, involving a dual approach via large deviation theory for occupation densities, can also be found in [13].

3.2. *The expected asymptotic shape.* We give a rough outline of calculations that will yield the correct exponential growth in the expected number of particles both near $-\gamma t$ in space and $\kappa\sqrt{t}$ in type at large times $t$. Using the formula from (17) and abusing notation throughout in the same way as Section 3.1, we find that

$$E \sum_{u \in N_t} \mathbb{1}_{\{X_u(t) \sim -\gamma t; Y_u(t) \geq \kappa\sqrt{t}\}}$$
$$= \mathbb{E}_{\mu_\lambda, \lambda}(e^{-\lambda \xi_t - \psi_\lambda^- \eta_t^2 + E_\lambda^- t} \mathbb{1}_{\{\xi_t \sim -\gamma t; \eta_t \geq \kappa\sqrt{t}\}})$$
$$\sim e^{(E_\lambda^- + \lambda\gamma - \kappa^2 \psi_\lambda^-)t} \mathbb{P}_{\mu_\lambda, \lambda}(\xi_t \sim -\gamma t; \eta_t \geq \kappa\sqrt{t}).$$

Now, from standard bounds on the tail of the normal distribution,

$$\mathbb{P}_{\mu_\lambda, \lambda}(\xi_t \sim -\gamma t; \eta_t \geq \kappa\sqrt{t})$$
(21)
$$= \mathbb{P}_{\mu_\lambda, \lambda}(\eta_t \geq \kappa\sqrt{t}) \mathbb{P}_{\mu_\lambda, \lambda}(\xi_t \sim -\gamma t | \eta_t \geq \kappa\sqrt{t})$$
$$\sim e^{-\mu_\lambda/\theta \kappa^2 t} \mathbb{P}_{\mu_\lambda, \lambda}(\xi_t \sim -\gamma t | \eta_t \geq \kappa\sqrt{t})$$

and, since $\psi_\lambda^- + (\mu_\lambda/\theta) = \psi_\lambda^+$, this yields

$$E \sum_{u \in N_t} \mathbb{1}_{\{X_u(t) \sim -\gamma t; Y_u(t) \geq \kappa\sqrt{t}\}}$$
(22)
$$\sim e^{(E_\lambda^- + \lambda\gamma - \kappa^2 \psi_\lambda^+)t} \mathbb{P}_{\mu_\lambda, \lambda}(\xi_t \sim -\gamma t | \eta_t \geq \kappa\sqrt{t}).$$

Recalling that $\Delta(\gamma, \kappa) := \inf_{\lambda \in (\lambda_{\min}, 0)} \{E_\lambda^- + \lambda\gamma - \kappa^2 \psi_\lambda^+\}$, simple calculus reveals this infimum is attained at a $\lambda$ value of

$$\bar{\lambda}(\gamma, \kappa) = -\gamma \sqrt{\frac{\theta(\theta - 8r)}{a^2(\kappa^2 + \theta)^2 + 4a\gamma^2\theta}} \in (\lambda_{\min}, 0), \tag{23}$$

and using this optimal value in equation (22) will lead to the upper bound

$$\limsup_{t \to \infty} t^{-1} \log E \sum_{u \in N_t} \mathbb{1}_{\{X_u(t) \sim -\gamma t; Y_u(t) \geq \kappa\sqrt{t}\}} \leq \Delta(\gamma, \kappa).$$

It is also clear from equation (22) that when minimizing $E_\lambda^- + \lambda\gamma - \kappa^2 \psi_\lambda^+$, we simultaneously maximize the probability $\mathbb{P}_{\mu_\lambda, \lambda}(\xi_t \sim -\gamma t | \eta_t \geq \kappa\sqrt{t})$. In particular, to get a matching lower bound, we do not want this probability to have any exponential decay in time when we choose the optimal parameter for $\lambda$.

In fact, at least up to the exponential decay rate in time, it can be shown using large-deviations arguments that

$$\mathbb{P}_{\mu_{\bar{\lambda}}, \bar{\lambda}}(\xi_t \sim -\gamma t; \eta_t \geq \kappa\sqrt{t}) \sim \exp\left(-\frac{\mu_{\bar{\lambda}}}{\theta}\kappa^2 t\right).$$

Indeed, we immediately gain the required upper bound from (21). For the lower bound, consider the following heuristics where we break paths into two sections: normal ergodic behavior over large time period $[0, t]$ followed by a rapid ascent out to type position $\kappa\sqrt{t}$ over a much shorter period $[t, t + \tau]$.

(i) *Ergodic behavior.* Over a large time $t$, the occupation density of $\eta$ will most likely have settled close to the invariant measure. Hence for large $t$, almost surely under $\mathbb{P}_{\lambda,\mu_\lambda}$,

$$\frac{1}{t}\int_0^t \eta_s^2\, ds \to \left(\frac{\theta}{2\mu_\lambda}\right).$$

(ii) *Rapid ascent.* Over a large time $\tau$, but where $\tau = o(t)$, the probability that $\eta$ starts close to the origin and ends near to $\kappa\sqrt{t}$, having followed close to the path $y$ over the entire time period $\tau$, is roughly given by

$$\exp\left(-\frac{1}{2\theta}\int_0^\tau \{\dot{y}(s) + \mu_\lambda y(s)\}^2\, ds\right)$$

under the $\mathbb{P}_{\lambda,\mu_\lambda}$ law. See, for example [24], Chapter 6, or [4], Chapter 5.6. After some Euler–Lagrange optimization, the path

$$y(s) = \kappa\sqrt{t}\frac{\sinh\mu_\lambda s}{\sinh\mu_\lambda \tau}$$

gives $\int_0^\tau y(s)^2\, ds \approx \kappa^2 t/(2\mu_\lambda)$, with the probability of this path being roughly $\exp(-(\mu_\lambda/\theta)\kappa^2 t)$.

Combining these two types of behavior, we can find paths with final positions

$$\eta_{t+\tau} \sim \kappa\sqrt{t}, \qquad \xi_{t+\tau} \sim \lambda a\int_0^{t+\tau} \eta_s^2\, ds \sim \lambda a\left(\frac{\theta}{2\mu_\lambda} + \frac{\kappa^2}{2\mu_\lambda}\right)t$$

and, moreover, when substituting the optimal $\lambda$ value of $\bar{\lambda}(\gamma, \kappa)$ and simplifying, this actually gives $\xi_{t+\tau} \sim -\gamma t$. Further, one of these paths occurs with a probability of roughly $\exp(-(\mu_{\bar{\lambda}}/\theta)\kappa^2 t)$, and note that $t + \tau \sim t$ since $\tau = o(t)$. Thus we see that to exponential order, the probability $\mathbb{P}_{\mu_{\bar{\lambda}},\bar{\lambda}}(\xi_t \sim -\gamma t; \eta_t \geq \kappa\sqrt{t})$ must be at least $\exp(-(\mu_{\bar{\lambda}}/\theta)\kappa^2 t)$, as required.

This heuristic argument can be made rigorous to prove, as claimed, that

$$\lim_{t\to\infty} t^{-1}\log E\left(\sum_{u\in N_t} \mathbb{1}_{\{X_u(t)\leq -\gamma t; Y_u(t)\geq \kappa\sqrt{t}\}}\right) = \Delta(\gamma, \kappa).$$

If we scale all spatial coordinates by $t^{-1}$ and all type coordinates by $(\sqrt{t})^{-1}$ at time $t$, the *expected asymptotic shape* can be considered to be the region $\mathcal{S} := \{(\gamma, \kappa) : \Delta(\gamma, \kappa) \geq 0\}$ where, on average, we have growth in the numbers of (scaled) particles.

**4. Short climb large deviation heuristics.** In this section, we give a heuristic calculation that suggests why the probability a single particle manages to have at least one descendant in the vicinity of $(-\beta t, \kappa\sqrt{t})$ near time $\tau$ is roughly $\exp(-\Theta(\beta,\kappa)t)$ for very large $t$, where $\Theta(\beta,\kappa)$ is given at equation (13). For these heuristics, we will think of $\tau$ as large and fixed, but of smaller order than $t$ (later on, in our rigorous approach, we will choose $\tau$ proportional to $\log t$). We emphasize that the heuristics in this section are neither meant to be precise nor made rigorous, yet they will provide invaluable intuition, guidance and motivation for our rigorous approach later on. Of particular importance will be the optimization problem that the heuristics suggest. Indeed, many of the exact calculations in Sections 4.2 and 4.3 will be essential later in the paper.

Suppose we start the branching diffusion with a single particle at $(0,0)$. First, we wish to know the probability that there is at least one particle at time $\tau$ that has a spatial position near $-\beta t$ having followed close to the path $x(s)$ for $0 \le s \le \tau$ *and* a type position near $\kappa\sqrt{t}$ having closely followed the path $y(s)$ for $0 \le s \le \tau$ for $t$ arbitrarily large.

We recall from large deviation theory of Ventcel–Freidlin (see [24], Chapter 6, or [4], Chapter 5.6) that the probability a *single* particle manages to follow closely both the type path $y(s)$ and the spatial path $x(s)$ for $0 \le s \le \tau$ is roughly given by

$$(24) \qquad \exp\left(-\frac{1}{2\theta}\int_0^\tau \left(\dot{y}(s) + \frac{\theta}{2}y(s)\right)^2 ds - \frac{1}{2}\int_0^\tau \frac{\dot{x}(s)^2}{ay(s)^2}\,ds\right)$$

when $x(0) = 0, x(\tau) = -\beta t$, $y(0) = 0, y(\tau) = \kappa\sqrt{t}$ and $t$ is very large. This probability will typically be very small, but if such paths are followed by particles in the branching diffusion, we have to also take account of the large breeding rates that are found far from the type origin.

If we let $X(s)$ represent the numbers of particles in the branching diffusion that are alive at time $s$ and have traveled "close" to the path $(x(u), y(u))$ for $0 \le u \le s$, then we can get a rough idea of how $X$ might behave by considering the following birth–death process.

4.1. *A birth–death process.* For *given* fixed paths $x(\cdot)$ and $y(\cdot)$, let $M$ be a time-dependent birth–death process where at time $s$ particles either give birth to single offspring with *breeding rate* $\lambda(s)$ given by

$$\lambda(s) = \rho + ry(s)^2,$$

or particles die with *death rate* $\mu(s)$ given by

$$\mu(s) = \frac{1}{2\theta}\left(\dot{y}(s) + \frac{\theta}{2}y(s)\right)^2 + \frac{1}{2}\frac{\dot{x}(s)^2}{ay(s)^2}.$$

(Note that the probability the initial particle of this birth–death process survives the entire time period $[0, \tau]$ is consistent with the rough large deviation probability for the branching diffusion at equation (24).)

An important quantity is the *effective total death rate* up to time $t$ which is defined by $\nu(s) := \int_0^s \{\mu(w) - \lambda(w)\} \, dw$, so here

$$\nu(s) = J(x, y, s)$$
$$:= \int_0^s \left( \frac{1}{2\theta} \left( \dot{y}(w) + \frac{\theta}{2} y(w) \right)^2 + \frac{1}{2} \frac{\dot{x}(w)^2}{ay(w)^2} - ry(w)^2 - \rho \right) ds.$$

The distribution for total number of offspring surviving, $M(\tau)$, for the time-dependent birth–death process is well known, for example, see [14]. Then defining

$$W_\tau := e^{-\nu(\tau)} \left( 1 + \int_0^\tau \mu(s) \, e^{\nu(s)} \, ds \right),$$
$$U_\tau := 1 - e^{-\nu(\tau)} W_\tau^{-1},$$
$$V_\tau := 1 - W_\tau^{-1},$$

we have

$$\mathbb{P}(M(\tau) = 0) = U_\tau,$$
$$\mathbb{P}(M(\tau) = n) = (1 - U_\tau)(1 - V_\tau) V_\tau^{n-1}, \qquad n = 1, 2, \ldots$$

with $\mathbb{E} M(\tau) = e^{-\nu(\tau)}$ and $\mathbb{E}(M(\tau) | M(\tau) \geq 1) = W_\tau$.

In our particular case, we have

$$\mathbb{E}(M(\tau)) = \exp(-J(x, y, \tau)).$$

Define the largest effective total death rate prior to time $\tau$ by

$$L(x, y, \tau) := \sup_{s \in [0, \tau]} J(x, y, s) \geq 0.$$

If we are in a case when $L(x, y, \tau)$ is *very* large, suggesting a high chance of extinction, then

$$(25) \qquad \mathbb{P}(M(\tau) \geq 1) = \frac{1}{1 + \int_0^\tau \mu(s) e^{\nu(s)} \, ds} \sim K_\tau \exp(-L(x, y, \tau)),$$

where $K_\tau^{-1} := \int_0^\tau \mu(s) \exp(-\{L(x, y, \tau) - J(x, y, s)\}) \, ds$. If there is at least one particle alive, we would then expect to have

$$\mathbb{E}(M(\tau) | M(\tau) \geq 1) \sim K_\tau^{-1} \exp(L(x, y, \tau) - J(x, y, \tau)).$$

Thus, we might *guess* that the probability any particles in the branching diffusion manage to make the difficult, rapid ascent along path $(x, y)$ to finish up near $(-\beta t, \kappa \sqrt{t})$ can, very *roughly*, be estimated by $\exp(-L(x, y, \tau))$. [To help see this, try writing $x(s) = tf(s)$ and $y(s) = \sqrt{t}g(s)$, thinking of $f, g$ as fixed paths and recall that $t$ is very large and $\tau = o(t)$, then the role of $K_\tau$ in (25) is insignificant next to $\exp(-L(x, y, \tau))$.]

We might then further *guess* that the chance any particles manage to stay near position $(-\beta t, \kappa\sqrt{t})$ during a very small interval of time close to $\tau$ should roughly look like

$$\exp\left(-\inf_{x,y} L(x, y, \tau)\right),$$

where we permit all possible paths $x$ and $y$ satisfying $x(0) = 0, x(\tau) = -\beta t$ and $y(0) = 0, y(\tau) = \kappa\sqrt{t}$ for the fixed time $\tau$. (We will state and prove a precise lower bound that corresponds to this guess at Theorem 7.)

4.2. *Finding the optimal path and probability.* We proceed to calculate

$$\inf_{x,y} L(x, y, \tau)$$

over paths $x$ and $y$ satisfying $x(0) = 0, x(\tau) = -\beta t$ and $y(0) = 0, y(\tau) = \kappa\sqrt{t}$ for the fixed time $\tau$.

We first note that

(26) $$\inf_{x,y} L(x, y, \tau) = \inf_{x,y} \sup_{s \in [0,\tau]} J(x, y, s) \geq \inf_{x,y} J(x, y, \tau)$$

and we now proceed to calculate $\inf_{x,y} J(x, y, \tau)$.

We can easily optimize over the choice of function $x$ given $y$, finding that

$$\dot{x}(s) \propto ay(s)^2 \quad \Rightarrow \quad x(s) = \lambda a \int_0^s y(u)^2 \, du$$

where $\lambda$ is the constant of proportionality and must satisfy

(27) $$\lambda = \frac{-\beta t}{a \int_0^\tau y(s)^2 \, ds},$$

yielding

$$\frac{1}{2} \int_0^\tau \frac{\dot{x}(s)^2}{ay(s)^2} \, ds = \frac{\beta^2 t^2}{2a \int_0^\tau y(s)^2 \, ds}.$$

This is exactly as anticipated since, when following the path $y$ in type space, the spatial position of a particle is following a Brownian motion with total amount of variance over period $\tau$ given by $a \int_0^\tau y(s)^2 \, ds$. Hence, the probability that a particle following the path $y$ in type space will also be found near to $\beta t$ in space at time $\tau$ is roughly

$$\exp\left(\frac{-\beta^2 t^2}{2a \int_0^\tau y(s)^2 \, ds}\right).$$

Introducing the notation
$$I(y) := \int_0^\tau \left( \frac{1}{2\theta} \left( \dot{y}(s) + \frac{\theta}{2} y(s) \right)^2 - r y(s)^2 \right) ds,$$
we are left to find

$$\inf_y \left\{ I(y) + \frac{\beta^2 t^2}{2a \int_0^\tau y(s)^2 \, ds} \right\}$$

(28)
$$= \inf_y \sup_\lambda \left\{ I(y) - \frac{1}{2} a \lambda^2 \int_0^\tau y(s)^2 \, ds - \lambda \beta t \right\}$$

$$\geq \sup_\lambda \inf_y \left\{ I(y) - \frac{1}{2} a \lambda^2 \int_0^\tau y(s)^2 \, ds - \lambda \beta t \right\}$$

where the first equality is trivially true by maximizing the quadratic in $\lambda$, the introduction of which conveniently removes the awkward integral in the denominator. Some further Euler–Lagrange optimization now gives the optimal path as

(29)
$$y_\lambda(s) = \kappa \sqrt{t} \frac{\sinh \mu_\lambda s}{\sinh \mu_\lambda \tau} \qquad (0 \leq s \leq \tau),$$

where
$$\mu_\lambda = \frac{\sqrt{\theta(\theta - 8r - 4a\lambda^2)}}{2},$$
and then

$$\sup_\lambda \inf_y \left\{ I(y) - \frac{1}{2} a \lambda^2 \int_0^\tau y(s)^2 \, ds - \lambda \beta t \right\}$$

$$= \sup_\lambda \left\{ \kappa^2 t \left( \frac{1}{4} + \frac{\mu_\lambda}{2\theta} \coth \mu_\lambda \tau \right) - \lambda \beta t \right\}.$$

The optimal parameter choice $\hat{\lambda}$ (which depends on $\tau$ as well as the model parameters) then satisfies

(30)
$$\frac{-\beta t}{a \hat{\lambda}} = \kappa^2 t \left( \frac{\coth \mu_{\hat{\lambda}} \tau}{2 \mu_{\hat{\lambda}}} - \frac{\tau}{2 \sinh^2 \mu_{\hat{\lambda}} \tau} \right) = \int_0^\tau y_{\hat{\lambda}}(s)^2 \, ds.$$

Then we have shown that
$$I(y_{\hat{\lambda}}) + \frac{\beta^2 t^2}{2a \int_0^\tau y_{\hat{\lambda}}(s)^2 \, ds} \geq \inf_y \left\{ I(y) + \frac{\beta^2 t^2}{2a \int_0^\tau y(s)^2 \, ds} \right\}$$

$$= \inf_y \sup_\lambda \left\{ I(y) - \frac{1}{2} a \lambda^2 \int_0^\tau y(s)^2 \, ds - \lambda \beta t \right\}$$

$$\geq \sup_\lambda \inf_y \left\{ I(y) - \frac{1}{2} a \lambda^2 \int_0^\tau y(s)^2 \, ds - \lambda \beta t \right\}$$

$$\geq I(y_{\hat{\lambda}}) - \frac{1}{2} a \hat{\lambda}^2 \int_0^\tau y_{\hat{\lambda}}(s)^2 \, ds - \hat{\lambda} \beta t,$$

and, in fact, we see that the left- and right-hand sides of the above are equal by recalling (30). It follows that the preceding supremum and infimum can be freely interchanged, actually preserving equality at the inequality (28).

Then, with the optimal spatial path

$$(31) \qquad x_\lambda(s) := \lambda a \int_0^s y_\lambda(u)^2 \, du = -\beta t \frac{\sinh 2\mu_\lambda s - 2\mu_\lambda s}{\sinh 2\mu_\lambda \tau - 2\mu_\lambda \tau},$$

and defining $\hat{x} := x_{\hat{\lambda}}, \hat{y} := y_{\hat{\lambda}}$, we have

$$\inf_{x,y} J(x,y,\tau) = J(\hat{x}, \hat{y}, \tau)$$

$$= t \sup_\lambda \left\{ \kappa^2 \left( \frac{1}{4} + \frac{\mu_\lambda}{2\theta} \coth \mu_\lambda \tau \right) - \lambda\beta \right\} - \rho\tau$$

$$= t \left( \kappa^2 \left( \frac{1}{4} + \frac{\mu_{\hat{\lambda}}}{2\theta} \coth \mu_{\hat{\lambda}} \tau \right) - \hat{\lambda}\beta \right) - \rho\tau.$$

Finally, it is easy to check that $J(\hat{x}, \hat{y}, \tau) = L(\hat{x}, \hat{y}, \tau)$, whence

$$\inf_{x,y} J(x,y,\tau) \geq \inf_{x,y} L(x,y,\tau),$$

and, combining with equation (26), we have found that

$$\inf_{x,y} L(x,y,\tau) = \inf_{x,y} J(x,y,\tau) = J(\hat{x}, \hat{y}, \tau).$$

4.3. *An important note on the optimal paths.* As $\tau \to \infty$, we have

$$\sup_\lambda \left\{ \kappa^2 \left( \frac{1}{4} + \frac{\mu_\lambda}{\theta} \coth \mu_\lambda \tau \right) - \lambda\beta \right\} \uparrow \sup_\lambda \{ \kappa^2 \psi_\lambda^+ - \lambda\beta \} = \kappa^2 \psi_{\bar{\lambda}}^+ - \bar{\lambda}\beta,$$

where the optimizing parameters of the supremums also converge with

$$(32) \qquad \hat{\lambda} \to \bar{\lambda} = -\beta \sqrt{\frac{\theta(\theta - 8r)}{a^2\kappa^4 + 4a\theta\beta^2}} = \bar{\lambda}\left( \left( \frac{\kappa^2 + \theta}{\kappa^2} \right) \beta, \kappa \right).$$

Note the agreement with previous optimal values at equations (23) and (15).

Then letting

$$(33) \qquad \Theta(\beta, \kappa) := \sup_\lambda \{ \kappa^2 \psi_\lambda^+ - \lambda\beta \}$$

$$= \frac{\kappa^2}{4} + \frac{\sqrt{\theta(\theta - 8r)(a^2\kappa^4 + 4a\theta\beta^2)}}{4a\theta}$$

and writing $\bar{x} := x_{\bar{\lambda}}$ and $\bar{y} := y_{\bar{\lambda}}$, we note that for all $\varepsilon, \delta > 0$ there exist $\tilde{\tau}, \mu > 0$ such that for all $t > 0$ and $\tau > \tilde{\tau}$

$$\exp\left( -\inf_{x,y} J(x,y,\tau) \right) \geq \exp(-J(\bar{x}, \bar{y}, \tau))$$

$$= \exp\left( -t \left( \kappa^2 \left( \frac{1}{4} + \frac{\mu_{\bar{\lambda}}}{2\theta} \coth \mu_{\bar{\lambda}} \tau \right) - \bar{\lambda}\beta \right) + \rho\tau \right)$$

$$\geq \exp(-t(\Theta(\beta,\kappa) + \varepsilon)).$$

Further (when $\kappa > 0$), for all $s \in [\tau - \mu, \tau]$,

$$\bar{y}(s) \geq (\kappa - \delta)\sqrt{t}, \qquad \bar{x}(s) \leq -(\beta - \delta)t.$$

In particular, the paths stay close to the required positions for some fixed length of time with corresponding probability at least as large as required.

**5. Proof of Theorem 3. Lower bound.** In this section we will state a precise short climb probability result and show how to combine it with almost sure spatial (only) growth rates to prove the lower bound of the growth rate in Theorem 3. This will make rigorous the two-phase mechanism described in Section 2 and suggested by the expectation calculations in Section 3.

The first phase requires knowledge of the almost-sure rates of growth of particles in the spatial dimension only. To this end, we will already make full use of Theorem 1 throughout this section, deferring its proof until Section 10.

The second phase requires a lower bound for the probability that a single particle makes a rapid ascent in type-space over the time interval $[0, \tau]$. This is the lower bound found in the heuristics of Section 4, but we require some further notation before the precise result can be stated. Note, throughout this section, we will only be interested in the optimal parameter value $\lambda = \bar{\lambda}$ as introduced in Section 4.3.

We wish to fix the relationship between *sufficiently large* $t$ and $\tau$ as

$$(34) \qquad \sqrt{\theta/(2\mu_{\bar{\lambda}})}\, e^{\mu_{\bar{\lambda}}\tau} = \kappa\sqrt{t}$$

and so define $\tau = \tau(t)$ by

$$(35) \qquad \tau(t) := \begin{cases} (2\mu_{\bar{\lambda}})^{-1} \log(2\mu_{\bar{\lambda}} t/\theta), & \text{for } 2\mu_{\bar{\lambda}} t > \theta, \\ 0, & \text{otherwise}. \end{cases}$$

Recall the optimal paths $(\bar{x}, \bar{y})$ over $s \in [0, \tau]$, where

$$(36) \qquad \bar{y}(s) = \kappa\sqrt{t}\, \frac{\sinh \mu_{\bar{\lambda}} s}{\sinh \mu_{\bar{\lambda}} \tau},$$

$$(37) \qquad \bar{x}(s) = a\bar{\lambda} \int_0^s \bar{y}(w)^2\, dw = -\beta t \frac{\sinh 2\mu_{\bar{\lambda}} s - 2\mu_{\bar{\lambda}} s}{\sinh 2\mu_{\bar{\lambda}} \tau - 2\mu_{\bar{\lambda}} \tau},$$

with fixed end points $\bar{y}(\tau) = \kappa\sqrt{t}$ and $\bar{x}(\tau) = -\beta t$.

For large times $t$ and $\delta, \varepsilon > 0$, let

$$(38) \quad A_t^{\varepsilon,\delta}(u) := \left\{ \sup_{s \in [0, \tau(t)]} |Y_u(s) - \bar{y}(s)| < \varepsilon\sqrt{t};\ \sup_{s \in [0, \tau(t)]} |X_u(s) - \bar{x}(s)| < \delta t \right\}.$$

We will use the notation

$$(39) \qquad A_t^{\varepsilon,\delta} := \bigcup_{u \in N_{\tau(t)}} A_t^{\varepsilon,\delta}(u)$$

for the event that there exists a particle in the branching diffusion that makes the short climb. Finally, recalling $\Theta(\beta, \kappa)$ given at (33), we can now state the short climb theorem:

THEOREM 7. *Fix any $y_1 > y_0 > 0$, $x \in \mathbb{R}$, and let $\varepsilon_0 > 0$. Then for any $\varepsilon$, $\delta > 0$, there exists $T > 0$ such that for all $y \in [y_0, y_1]$,*

$$t^{-1} \log P^{x,y}(A_t^{\varepsilon,\delta}) \geq -(\Theta(\beta, \kappa) + \varepsilon_0)$$

*for all $t > T$.*

We will prove Theorem 7 using a spine change of measure. This requires us to introduce the notation for the spine set-up in detail before proceeding, so this and further technical issues are postponed to Sections 6 and 7.

REMARK 8. We note that Theorem 7 is actually a stronger result than needed to prove Theorem 3 because we identify the specific *paths* followed by particles that are near position $(\beta t, \kappa \sqrt{t})$ at time $t + \tau$, rather than just considering the particle's *positions* close to time $t + \tau$.

In combining the two phases, we will have a huge number of independent trials each with a small probability of success, intuitively giving rise to a Poisson approximation for a large number of successful particles. In fact, in our proof of the lower bound of Theorem 3 below, we will actually use the following result about the behavior of sequences of sums of independent Bernoulli random variables.

LEMMA 9. *For each $n$, define the random variable $B_n := \sum_{u \in F_n} \mathbb{1}_{E_n(u)}$ where the events $\{E_n(u) : u \in F_n\}$ are independent. Let $p_n(u) := P(E_n(u))$ and $S_n := \sum_{u \in F_n} p_n(u)$ and suppose that, for some $\nu \in (1/2, 1)$,*

$$\sum_{n \in \mathbb{N}} \frac{1}{(S_n)^{2\nu - 1}} < \infty. \tag{40}$$

*Then the sequence of (possibly dependent) random variables $\{B_1, B_2, \ldots\}$ has $|B_n - S_n| > (S_n)^\nu$ for only finitely many $n$, almost surely.*

*In particular, for any $\varepsilon > 0$, there exists some (random) $N \in \mathbb{N}$ such that, with probability one,*

$$\frac{B_n}{S_n} > 1 - \varepsilon \qquad \text{for all } n > N. \tag{41}$$

PROOF. For $\nu \in (1/2, 1)$, Chebyshev's inequality yields

$$\mathbb{P}(|B_n - S_n| > (S_n)^\nu) \leq \frac{\sum_{u \in F_n} p_n(u)(1 - p_n(u))}{(S_n)^{2\nu}} \leq \frac{1}{(S_n)^{2\nu - 1}},$$

and hence the Borel–Cantelli lemmas, combined with hypothesis (40), imply that

$$|B_n - S_n| > (S_n)^\nu$$

for only finitely many $n$, almost surely. Equation (41) now follows on division by $S_n$, and noticing the assumption (40) implies that $\lim_{n\to\infty} S_n = \infty$. □

PROOF OF THEOREM 3. LOWER BOUND. Define $f^{-1}(t) := t - \tau(t)$, noting that both $f(t)/t \to 1$ and $f^{-1}(t)/t \to 1$ as $t \to \infty$. Also, for $n \in \mathbb{N}$ and $\mu > 0$, define $T_n := (n+1)\mu$. We want to estimate the number of particles that are near the large position $(-(\alpha + \beta)T_n, \kappa\sqrt{T_n})$ during time interval $[T_{n-1}, T_n]$. For this, we will consider particles that travel with a velocity $-\alpha$ over time period $[0, f^{-1}(T_n)]$ before commencing their rapid ascent of (relatively short) duration $\tau(T_n)$ to be in final position at time $T_n$. Then

$$\inf_{s \in [T_{n-1}, T_n]} N_s((\alpha + \beta - \delta)T_n; [(\kappa - \delta)\sqrt{T_n}, \infty))$$

$$(42) \qquad \geq \sum_{u \in N_{T_n}} \mathbb{1}_{\{\bigcap_{s \in [T_{n-1}, T_n]}\{X_u(s) \leq -(\alpha+\beta-\delta)T_n; Y_u(s) \geq (\kappa-\delta)\sqrt{T_n}\}\}}$$

$$\geq \sum_{u \in F_n^\alpha} \mathbb{1}_{\{\bar{N}_n^{\beta,\kappa}(u) > 0\}}$$

where

$$F_n^\alpha := \{u \in N_{f^{-1}(T_n)} : X_u(f^{-1}(T_n)) \leq -\alpha T_n, Y_u(f^{-1}(T_n)) \in [y_0, y_1]\}$$

and, for $u \in F_n^\alpha$,

$$\bar{N}_n^{\beta,\kappa}(u) := \sum_{\substack{v \in N_{T_n} \\ v \geq u}} \mathbb{1}_{\{\bigcap_{s \in [T_{n-1}, T_n]}\{X_v(s) - X_v(f^{-1}(T_n)) \leq -(\beta-\delta)T_n; Y_v(s) \geq (\kappa-\delta)\sqrt{T_n}\}\}}.$$

We will now show that the sum at (42) grows as fast as anticipated:

LEMMA 10. *For any $\varepsilon > 0$, we may choose $\mu > 0$ such that there exists a random $N \in \mathbb{N}$ where*

$$\frac{1}{T_n} \log \sum_{u \in F_n^\alpha} \mathbb{1}_{\{\bar{N}_n^{\beta,\kappa}(u) > 0\}} \geq \Delta(\alpha) - \Theta(\beta, \kappa) - \varepsilon$$

*for all $n > N$ with probability one.*

PROOF. We will be able to apply Lemma 9 given sufficient information about the growth of $|F_n^\alpha|$ and decay of the probabilities

$$p_n^{\beta,\kappa}(u) := P(\bar{N}_n^{\beta,\kappa}(u) > 0 | \mathcal{F}_{f^{-1}(T_n)}),$$

where $u \in F_n^\alpha \subset N_{f^{-1}(T_n)}$.

It follows easily from Theorem 1, $f^{-1}(T_n)/T_n \to 1$ and the continuity of $\Delta(\alpha)$ that

$$\frac{\log |F_n^\alpha|}{T_n} \geq \Delta(\alpha) - \frac{\varepsilon}{4}$$

for all sufficiently large $n$.

The definition of $\bar{N}_n^{\beta,\kappa}(u)$ and spatial translation invariance implies that, for each $u \in F_n^\alpha$, the rapid ascent probability $p_n^{\beta,\kappa}(u)$ depends only on the initial type position $Y_u(f^{-1}(T_n))$.

For $\delta, \mu > 0$, define

$$B_t^{\delta,\mu}(u) := \bigcap_{s \in [\tau(t)-\mu,\tau(t)]} \{X_u(s) - X_u(0) < -(\beta-\delta)t; Y_u(s) \geq (\kappa-\delta)\sqrt{t}\}$$

and

(43) $$B_t^{\delta,\mu} := \bigcup_{u \in N_{\tau(t)}} B_t^{\delta,\mu}(u).$$

Recalling the comments of Section 4.3, there exist $\varepsilon', \delta' > 0$ and we may choose $\mu > 0$ sufficiently small, such that

$$p_n^{\beta,\kappa}(u) = P^{0,Y_u(f^{-1}(T_n))}(B_{T_n}^{\delta,\mu}) \geq P^{0,Y_u(f^{-1}(T_n))}(A_{T_n}^{\varepsilon',\delta'}) =: \bar{p}_n(u)$$

for all $u \in F_n^\alpha$ whenever $n$ is sufficiently large. Together with Theorem 7 and since $Y_u(f^{-1}(T_n)) \in [y_0, y_1]$ for $u \in F_n^\alpha$, this reveals

$$\frac{\log p_n^{\beta,\kappa}(u)}{T_n} \geq \frac{\log \bar{p}_n(u)}{T_n} \geq -\Theta(\beta,\kappa) - \frac{\varepsilon}{4}$$

for all for $u \in F_n^\alpha$ and all sufficiently large $n$, almost surely. Then we may combine the observations above to obtain

$$\frac{1}{T_n} \log \sum_{u \in F_n^\alpha} p_n^{\beta,\kappa}(u) \geq \Delta(\alpha) - \Theta(\beta,\kappa) - \frac{\varepsilon}{2}.$$

Taking this last line together the assertion of Lemma 9 at equation (41) gives the result. $\square$

It is now straightforward to combine Lemma 10 with the inequality at (42) to see that, given $\varepsilon, \delta > 0$, there exists $\mu > 0$ and a random time $T$ such that

$$t^{-1} \log N_t((\alpha+\beta-\delta)t; [(\kappa-\delta)\sqrt{t}, \infty)) \geq \Delta(\alpha) - \Theta(\beta,\kappa) - \varepsilon$$

for all $t > T$, almost surely. Since $\varepsilon$ and $\delta$ can be taken arbitrarily small, using the optimal $\bar{\alpha}$ and $\bar{\beta}$ according to equations (14)–(15), we find

$$\liminf_{t \to \infty} t^{-1} \log N_t(\gamma, [\kappa\sqrt{t}, \infty)) \geq \Delta(\gamma, \kappa) \quad \text{almost surely,}$$

as required. [It is also interesting to note that $\bar{\lambda} = \lambda_{\bar{\alpha}} = \bar{\lambda}(\gamma, \kappa)$ from equations (19), (23) and (32), so the optimal parameters are in agreement with those of the expectation calculations in Section 3 and the path large deviations in Section 4.] □

**6. The "spine" setup and results.** In this section, we describe how to construct an enriched branching diffusion with an identified "spine" or "backbone" particle and discuss how to perform some extremely useful changes of measure (closely related to the additive martingales) that will essentially "force" the spine perform the short climb, whilst giving birth at an accelerated rate to offspring that behave as if under the original measure. These spine techniques are at the very heart of our proof of Theorem 7 in Section 5. Spine ideas were first seen for branching Brownian motion in [3] and developed for Galton–Watson processes in [16, 18, 19]. Kyprianou [17] and Englander and Kyprianou [6], developed the technique for some families of branching diffusions; and more recently the spine approach has been significantly improved in [8]. This approach uses several different filtrations on an enlarged probability space carrying the branching diffusion, and permits some very useful techniques and results to be developed. For example, "additive" (many-particle) martingales can be represented as suitable conditional expectations of "spine" (single-particle) martingales and consequently there are clear interpretations for any changes of measure and *all* measures involved in our "spine" setup are *probability* measures with intuitive constructions. Following Hardy and Harris [8], we will first outline the notation and then describe the changes of measure. The notation described in this section is generalized to allow each particle $u$ to have $1 + A_u$ offspring, where each $A_u$ is an independent copy of a random variable with values in $\{0, 1, 2, \ldots\}$. The spine techniques developed in this paper could readily be generalized to such models.

All probability measures are to be defined on the space $\tilde{\mathcal{T}}$ of *marked Galton–Watson trees with spines*; before defining precisely what this space is we need to set up some other notation. We recall the set of Ulam–Harris labels, $\Omega$, defined by $\Omega := \{\varnothing\} \cup \bigcup_{n \in \mathbb{N}} (\mathbb{N})^n$, where $\mathbb{N} := \{1, 2, 3, \ldots\}$. For two words $u, v \in \Omega$, $uv$ denotes the concatenated word, where we take $u\varnothing = \varnothing u = u$. So $\Omega$ contains elements such as "$\varnothing 412$," which represents "the individual being the 2nd child of the 1st child of the 4th child of the initial ancestor $\varnothing$." For labels $u, v \in \Omega$ the notation $v < u$ means that $v$ is an ancestor of $u$, and $|u|$ denotes the length of $u$.

We define a *Galton–Watson tree* to be a set $\tau \subset \Omega$ such that:

(i) $\varnothing \in \tau$, so there is the unique initial ancestor;

(ii) if $u, v \in \Omega$, then $vu \in \tau \Rightarrow v \in \tau$, so $\tau$ contains all of the ancestors of its nodes;

(iii) for all $u \in \tau$, there exists $A_u \in \{0, 1, 2, \ldots\}$ such that for $j \in \mathbb{N}$, $uj \in \tau$ if and only if $1 \leq j \leq 1 + A_u$.

The set of all such trees is $\mathbb{T}$, and we will use the symbol $\tau$ for a particular tree. As our work concerns branching diffusions we shall often refer to the labels of $\tau$ as particles. Note that for the *binary* branching mechanism in this paper, $P(A_u = 1) \equiv 1$; of course, here there is only one $\tau \in \mathbb{T}$—the binary tree.

A Galton–Watson tree by itself only records the family structure of the individuals, so to each individual $u \in \tau$ we give a mark $(X_u, Y_u, \sigma_u)$ which contains the following information:

- $\sigma_u \in [0, \infty)$ is the *lifetime* of particle $u$, which also determines the *fission time* of the particle as $S_u := \sum_{v \leq u} \sigma_v$. We may also refer to the $S_u$ as *death times*;
- the function $X_u(t) : [S_u - \sigma_u, S_u) \to \mathbb{R}$ describes the particle's *spatial* motion in $\mathbb{R}$ during its lifetime;
- the function $Y_u(t) : [S_u - \sigma_u, S_u) \to \mathbb{R}$ describes the evolution of the particle's *type* in $\mathbb{R}$ during its lifetime.

For clarity we must decide whether or not a particle is in existence at its death time: our convention will be that a particle dies "infinitesimally before" its death time—this is why $X_u$ and $Y_u$ are defined on $[S_u - \sigma_u, S_u)$ and not $[S_u - \sigma_u, S_u]$—so that at time $S_u$ the particle $u$ has disappeared and has been replaced by its two children.

We denote a particular marked tree by $(\tau, X, Y, \sigma)$, or the abbreviation $(\tau, M)$, and the set of all marked Galton–Watson trees by $\mathcal{T}$. For each $(\tau, X, Y, \sigma) \in \mathcal{T}$, the set of particles alive at time $t$ is defined as $N_t := \{u \in \tau : S_u - \sigma_u \leq t < S_u\}$. For any given marked tree $(\tau, M) \in \mathcal{T}$ we can distinguish individual lines of descent from the initial ancestor: $\varnothing, u_1, u_2, u_3, \ldots \in \tau$, where $u_i$ is a child of $u_{i-1}$ for all $i \in \{2, 3, \ldots\}$ and $u_1$ is a child of the initial individual $\varnothing$. We call such a line of descent a *spine* and denote it by $\xi$. In a slight abuse of notation we refer to $\xi_t$ as the unique node in $\xi$ that is alive at time $t$, and also for the position of the particle that makes up the spine at time $t$; that is, $\xi_t := X_u(t)$, where $u \in \xi \cap N_t$. However, although the interpretation of $\xi_t$ should always be clear from the context, we introduce the following notation for use where some ambiguity may still arise:

- $\mathrm{node}_t((\tau, M, \xi)) := u$ if $u \in \xi$ is the node in the spine alive at time $t$.

It is natural to think of the spine as a single diffusing particle $\xi_t$, or, strictly speaking, the pair $(\xi_t, \eta_t)$, where $\eta_t$ is the type of the spine at time $t$.

We define $n_t$ to be a counting function that tells us which generation of the spine is currently alive, or equivalently the number of fission times there have been on the spine:

$$n_t = |\mathrm{node}_t(\xi)|.$$

The collection of all marked trees with a distinguished spine is the space $\tilde{\mathcal{T}}$ on which our probability measures will eventually be defined, but first we define four filtrations on this space that contain different levels of information about the branching diffusion.

- Filtration $(\mathcal{F}_t)_{t\geq 0}$. We define a filtration of $\tilde{\mathcal{T}}$ made up of the $\sigma$-algebras

$$\mathcal{F}_t := \sigma((u, X_u, Y_u, \sigma_u) : S_u \leq t;$$
$$(u, X_u(s), Y_u(s) : s \in [S_u - \sigma_u, t]) : t \in [S_u - \sigma_u, S_u)),$$

which means that $\mathcal{F}_t$ is generated by the information concerning all particles that have lived and died before time $t$, and also those that are still alive at time $t$. Each of these $\sigma$-algebras is a subset of the limit $\mathcal{F}_\infty := \sigma(\bigcup_{t\geq 0} \mathcal{F}_t)$.

- Filtration $(\tilde{\mathcal{F}}_t)_{t\geq 0}$. We define the filtration $(\tilde{\mathcal{F}}_t)_{t\geq 0}$ by augmenting the filtration $\mathcal{F}_t$ with the knowledge of which node is the spine at time $t$; that is, $(\tilde{\mathcal{F}}_t)_{t\geq 0} := \sigma(\mathcal{F}_t, \text{node}_t(\xi))$ and $\tilde{\mathcal{F}}_\infty := \sigma(\bigcup_{t\geq 0} \tilde{\mathcal{F}}_t)$, so that this filtration knows everything about the branching diffusion and everything about the spine.

- Filtration $(\mathcal{G}_t)_{t\geq 0}$. $(\mathcal{G}_t)_{t\geq 0}$ is a filtration of $\tilde{\mathcal{T}}$ defined by $\mathcal{G}_t := \sigma(\xi_s : 0 \leq s \leq t)$, and $\mathcal{G}_\infty := \sigma(\bigcup_{t\geq 0} \mathcal{G}_t)$. These $\sigma$-algebras are generated only by the spine's motion and so do not contain the information about which nodes of the tree $\tau$ make up the spine.

- Filtration $(\tilde{\mathcal{G}}_t)_{t\geq 0}$. As we did in going from $\mathcal{F}_t$ to $\tilde{\mathcal{F}}_t$ we create $(\tilde{\mathcal{G}}_t)_{t\geq 0}$ from $(\mathcal{G}_t)_{t\geq 0}$ by including knowledge of which nodes make up the spine: $(\tilde{\mathcal{G}}_t)_{t\geq 0} := \sigma(\mathcal{G}_t, \text{node}_t(\xi))$ and $\tilde{\mathcal{G}}_\infty := \sigma(\bigcup_{t\geq 0} \tilde{\mathcal{G}}_t)$. This means that $\tilde{\mathcal{G}}_t$ also knows when the fission times on the spine occurred, whereas $\mathcal{G}_t$ does not.

Now that we have defined the underlying space and filtrations, we can define the probability measures of interest. We let the typed branching diffusion be as described in Section 1.1, with the probability measures $\{P^{x,y} : x, y \in \mathbb{R}\}$ on $(\tilde{\mathcal{T}}, \mathcal{F}_\infty)$ representing the law of this typed branching diffusion when initially started with a single particle at $(x, y)$.

We recall from [18] that, if $f$ is an $\tilde{\mathcal{F}}_t$-measurable function, we can write

(44) $$f = \sum_{u \in N_t} f_u \mathbb{1}_{\{\xi_t = u\}},$$

where $f_u$ is $\mathcal{F}_t$-measurable. Now we can extend $P^{x,y}$ to a measure $\tilde{P}^{x,y}$ on $(\tilde{\mathcal{T}}, \tilde{\mathcal{F}}_\infty)$ by choosing the particle that continues the spine uniformly each time there is a birth on the spine; more precisely, for any $f \in m\tilde{\mathcal{F}}_t$ with representation like (44), we have:

$$\int_{\tilde{\mathcal{T}}} f \, d\tilde{P}^{x,y}(\tau, M, \xi) := \int_{\mathcal{T}} \sum_{u \in N_t} f_u \prod_{v < u} \tfrac{1}{2} \, dP^{x,y}(\tau, M).$$

We construct the $\tilde{\mathcal{F}}_t$-measurable martingale $\tilde{\zeta}(t)$ as

$$\tilde{\zeta}(t) := v_\lambda^+(\eta_t) e^{\int_0^t \{R(\eta_s) + 1/2\lambda^2 A(\eta_s)\} \, ds - E_\lambda^+ t} \times 2^{n_t} e^{-\int_0^t R(\eta_s) \, ds}$$

$$(45) \qquad \times e^{\lambda \xi_t - 1/2\lambda^2 \int_0^t A(\eta_s) \, ds}$$

$$= v_\lambda^+(\eta_t) 2^{n_t} e^{\lambda \xi_t - E_\lambda^+ t}.$$

Observe that this is a product of single-particle martingales, details of which can be found in [17] or [10]. One can think of these as $h$-transforms of the $\tilde{P}$-law of the spine: the first makes $\eta$ an outward-drifting Ornstein–Uhlenbeck process with drift parameter $\mu_\lambda$; the second increases the breeding rate on the spine to $2R(\cdot)$; and the third adds a spatial drift to $\xi$.

Using the martingale $\tilde{\zeta}(t)$ we may define a measure $\tilde{\mathbb{Q}}_\lambda^{x,y}$ on $(\tilde{\mathcal{T}}, \tilde{\mathcal{F}}_\infty)$ by

$$(46) \qquad \left. \frac{d\tilde{\mathbb{Q}}_\lambda^{x,y}}{d\tilde{P}^{x,y}} \right|_{\tilde{\mathcal{F}}_t} = \frac{\tilde{\zeta}(t)}{\tilde{\zeta}(0)} = \frac{e^{-\lambda x}}{v_\lambda^+(y)} v_\lambda^+(\eta_t) 2^{n_t} e^{\lambda \xi_t - E_\lambda^+ t}.$$

And since $\tilde{\zeta}(t)$ is a product of $h$-transforms, under $\tilde{\mathbb{Q}}_\lambda^{x,y}$ the process may be re-constructed path-wise according to the following description:

- starting from spatial position $x$ and type $y$ the spine $(\xi_t, \eta_t)$ diffuses spatially as a Brownian motion with infinitesimal variance $A(\eta_t)$ and infinitesimal drift $\lambda A(\eta_t)$;
- the type of the spine, $\eta_t$, begins at $y$ and moves in type space as an *outward-drifting* Ornstein–Uhlenbeck process with generator

$$\frac{\theta}{2} \frac{\partial^2}{\partial y^2} + \mu_\lambda y \frac{\partial}{\partial y};$$

- the spine branches at rate $2R(\eta_t)$, producing 2 particles;
- one of these particles is selected uniformly at random;
- the chosen offspring repeats stochastically the behavior of its parent;
- the other offspring particle initiates a $P^{\cdot,\cdot}$-BBM from its birth position and type.

The change of measure (46) projects onto the sub-algebra $\mathcal{F}_t$ as a conditional expectation:

$$\left. \frac{d\tilde{\mathbb{Q}}_\lambda^{x,y}}{d\tilde{P}^{x,y}} \right|_{\mathcal{F}_t} = \frac{e^{-\lambda x}}{v_\lambda^+(y)} \tilde{P}^{x,y}(v_\lambda^+(\eta_t) 2^{n_t} e^{\lambda \xi_t - E_\lambda^+ t} | \mathcal{F}_t),$$

and it is a short calculation using the methods of, for example, Hardy and Harris [10] to show that:

THEOREM 11. *If we define $\mathbb{Q}_\lambda^{x,y} := \tilde{\mathbb{Q}}_\lambda^{x,y}|_{\mathcal{F}_\infty}$, then $\mathbb{Q}_\lambda^{x,y}$ is a measure on $\mathcal{F}_\infty$ that satisfies*

$$\left.\frac{d\mathbb{Q}_\lambda^{x,y}}{dP^{x,y}}\right|_{\mathcal{F}_t} = \hat{Z}_\lambda^+(t) := \frac{Z_\lambda^+(t)}{Z_\lambda^+(0)}.$$

*Moreover under $\mathbb{Q}_\lambda^{x,y}$, the path-wise construction of the branching diffusion is the same as under $\tilde{\mathbb{Q}}_\lambda$.*

Although the path-wise construction of the branching diffusion is the same under $\mathbb{Q}_\lambda^{x,y}$ and $\tilde{\mathbb{Q}}_\lambda^{x,y}$, only the measure $\tilde{\mathbb{Q}}_\lambda^{x,y}$ "knows" about the spine. It is clear, however, that we have $\tilde{\mathbb{Q}}_\lambda^{x,y}(A) = \mathbb{Q}_\lambda^{x,y}(A)$ for any $A \in \mathcal{F}_\infty$.

Under the measure $\tilde{\mathbb{Q}}_\lambda^{x,y}$ only the behavior of the spine is altered, and combining this observation with conditioning on the spine's path and fission-times gives us a very useful representation for $Z_\lambda^+(t)$ under $\tilde{\mathbb{Q}}_\lambda^{x,y}$ that we shall refer to as the *spine decomposition*:

$$(47) \quad \tilde{\mathbb{Q}}_\lambda^{x,y}(Z_\lambda^+(t)|\tilde{\mathcal{G}}_\infty) = \sum_{u<\xi_t} v_\lambda^+(\eta_{S_u})e^{\lambda\xi_{S_u} - E_\lambda^+ S_u} + v_\lambda^+(\eta_t)e^{\lambda\xi_t - E_\lambda^+ t}.$$

Throughout the rest of this article we will refer to the two pieces of this decomposition as the "sum term" and the "spine term." This decomposition is discussed in detail for a wide variety of branching diffusions in [9], but to derive it we simply note that the contributions to $Z_\lambda^+(t)$ from the subtrees that branch off the spine have *constant* $\tilde{\mathbb{Q}}_\lambda^{x,y}$-expectation because they behave as if under the original measure $P$, and we know that $Z_\lambda^+(t)$ is a $P$-martingale. The spine decomposition reduces many calculations about the behavior of $Z_\lambda^+(t)$ under $\tilde{\mathbb{Q}}_\lambda^{x,y}$ to one-particle calculations about the spine, and this observation is exploited in the spine proofs of $\mathcal{L}^p$-bounds for some families of additive martingales in [9].

**7. Proof of Theorem 7. The short climb probability.** With the spine foundations firmly established in Section 6, we may proceed with the proof of the short climb probability lower bound from Theorem 7.

First, recall definitions (38) and (39), where $A_t^{\varepsilon,\delta}$ is the event that *there exists* a particle that makes the short climb along optimal path $(\bar{x}, \bar{y})$, and $A_t^{\varepsilon,\delta}(\xi)$ is the event that the *spine* makes the short climb. Note that $\varepsilon$ controls the proximity to $\bar{x}$ and $\delta$ the proximity to $\bar{y}$. Importantly, we will *only* be interested in taking $\lambda = \bar{\lambda}$ throughout this section, although we will usually just write $\lambda$ for notational simplicity. Also recall throughout that $t$ and $\tau$ are related through $(\theta/(2\mu_\lambda))\exp(2\mu_\lambda \tau) = \kappa^2 t$.

PROOF OF THEOREM 7. The key step in the proof of this is the following use of the spine change of measure: for any function $g:\mathbb{R}^+ \to \mathbb{R}^+$ we have

$$P^{x,y}(A_t^{\varepsilon,\delta}) = \mathbb{Q}_\lambda^{x,y}\left(\frac{\mathbb{1}_{A_t^{\varepsilon,\delta}}}{\hat{Z}_\lambda^+(\tau)}\right) = \tilde{\mathbb{Q}}_\lambda^{x,y}\left(\frac{\mathbb{1}_{A_t^{\varepsilon,\delta}}}{\hat{Z}_\lambda^+(\tau)}\right) \geq \tilde{\mathbb{Q}}_\lambda^{x,y}\left(\frac{\mathbb{1}_{A_t^{\varepsilon,\delta}(\xi)}}{\hat{Z}_\lambda^+(\tau)}\right)$$

(48)
$$\geq \tilde{\mathbb{Q}}_\lambda^{x,y}\left(\frac{\mathbb{1}_{A_t^{\varepsilon,\delta}(\xi)}}{\hat{Z}_\lambda^+(\tau)}; \sup_{s\in[0,\tau]} \hat{Z}_\lambda^+(s) \leq g(\tau)\right)$$

$$\geq g(\tau)^{-1}\tilde{\mathbb{Q}}_\lambda^{x,y}\left(A_t^{\varepsilon,\delta}(\xi); \sup_{s\in[0,\tau]} \hat{Z}_\lambda^+(s) \leq g(\tau)\right).$$

Essentially we just have to make the "correct" choice for both $\lambda$ and $g$ in expression (48), although there will still remain a number of technicalities to resolve.

The first idea is to ensure the (originally rare) event $A_t^{\varepsilon,\delta}$ actually occurs under the new measure $\tilde{\mathbb{Q}}_\lambda^{x,y}$ by making the spine follow close to the required path $(\bar{x}, \bar{y})$; this is achieved by choosing the optimal value $\bar{\lambda}$ for $\lambda$ and choosing $\tau$ to be on the natural time scale it would take the spine to reach position $\kappa\sqrt{t}$. In particular, this choice will mean that in the first line of the above set of inequalities there is no significant loss of mass when replacing the event $A_t^{\varepsilon,\delta}$ with $A_t^{\varepsilon,\delta}(\xi)$. Next, we wish to choose the smallest possible $g$ that will still leave some positive probability on the last line of the above argument. Hence, we wish to identify the rate of growth of the martingale $Z_\lambda^+$ under $\tilde{\mathbb{Q}}_\lambda^{x,y}$, and this will essentially be governed by the contribution from the spine itself.

With this is mind, and recalling the various properties of the optimal paths and parameters from Section 4, for $\varepsilon_0 > 0$ we define

$$g_{\varepsilon_0}(\tau) := \exp\left(\left(\psi_\lambda^+ + \frac{\lambda^2 a}{2\mu_\lambda} + \frac{\varepsilon_0}{\kappa}\right)\left(\frac{\theta}{2\mu_\lambda}\right)e^{2\mu_\lambda \tau} - (\psi_\lambda^+ y^2 + \lambda x)\right)$$

and recall from (35) that the scaling between $t$ and $\tau$ is fixed throughout, where $\kappa^2 t = (\theta/(2\mu_{\bar{\lambda}}))e^{2\mu_{\bar{\lambda}}\tau}$ for large $t$, hence $t + \tau \sim t$. Note that since we are only considering the optimal value $\lambda = \bar{\lambda}$, we have

$$\left(\psi_{\bar{\lambda}}^+ + \frac{\bar{\lambda}^2 a}{2\mu_{\bar{\lambda}}}\right)\left(\frac{\theta}{2\mu_{\bar{\lambda}}}\right)e^{2\mu_{\bar{\lambda}}\tau} = (\kappa^2\psi_{\bar{\lambda}}^+ - \bar{\lambda}\beta)t = \Theta(\beta, \kappa).$$

Then from (48) we have

(49) $$P^{x,y}(A_t^{\varepsilon,\delta}) \geq g_{\varepsilon_0}(\tau)^{-1}\tilde{\mathbb{Q}}_\lambda^{x,y}\left(A_t^{\varepsilon,\delta}(\xi); \sup_{s\in[0,\tau]} \hat{Z}_\lambda^+(s) \leq g_{\varepsilon_0}(\tau)\right).$$

Our strategy for the rest of this proof is to show that the $\tilde{\mathbb{Q}}_\lambda^{x,y}$-probability in (49) is at least some $\varepsilon' > 0$ for all sufficiently large $t$, uniformly for $y \in [y_0, y_1]$, so that the decay rate part of (49) matches the desired rate in the statement of the theorem.

Conditioning on the spine's path and birth times, $\tilde{\mathcal{G}}_\infty$, and then making use of some standard properties of conditional expectation we have

$$\tilde{\mathbb{Q}}_\lambda^{x,y}\left(A_t^{\varepsilon,\delta}(\xi); \sup_{s\in[0,\tau]} \hat{Z}_\lambda^+(s) \leq g_{\varepsilon_0}(\tau)\right)$$

$$= \tilde{\mathbb{Q}}_\lambda^{x,y}\left(\tilde{\mathbb{Q}}_\lambda^{x,y}\left(A_t^{\varepsilon,\delta}(\xi); \sup_{s\in[0,\tau]} \hat{Z}_\lambda^+(s) \leq g_{\varepsilon_0}(\tau)\Big|\tilde{\mathcal{G}}_\infty\right)\right)$$

$$= \tilde{\mathbb{Q}}_\lambda^{x,y}\left(\mathbb{1}_{A_t^{\varepsilon,\delta}(\xi)} \tilde{\mathbb{Q}}_\lambda^{x,y}\left(\sup_{s\in[0,\tau]} \hat{Z}_\lambda^+(s) \leq g_{\varepsilon_0}(\tau)\Big|\tilde{\mathcal{G}}_\infty\right)\right),$$

since $A_t^{\varepsilon,\delta}(\xi)$ is $\tilde{\mathcal{G}}_\infty$-measurable. We next observe that, conditional on $\tilde{\mathcal{G}}_\infty$, we can write $\hat{Z}_\lambda^+(t)$ as

$$(50) \qquad \hat{Z}_\lambda^+(t) = e^{-(\psi_\lambda^+ y^2+\lambda x)}\left(\sum_{u<\xi_t} e^{-E_\lambda^+ S_u} Z_\lambda^{(u)}(t-S_u) + f(t)\right),$$

where the $Z_\lambda^{(u)}$ are independent copies of $Z_\lambda^+$ started from a single particle at $(\xi_{S_u}, \eta_{S_u})$; and $f(t)$ is the contribution to $Z_\lambda^+(t)$ from the spine, which, conditional on $\tilde{\mathcal{G}}_\infty$, is a known function of $t$. Now if we could show, for $0 < \tilde{\varepsilon}_0 < \varepsilon_0$,

$$\sup_{s\in[0,\tau]} \hat{f}(s) \leq \frac{g_{\tilde{\varepsilon}_0}(\tau)}{2} \quad \text{and} \quad \sup_{s\in[0,\tau]} (\hat{Z}_\lambda^+(s) - \hat{f}(s)) \leq \frac{g_{\varepsilon_0}(\tau)}{2},$$

where $\hat{f}(t) := e^{-(\psi_\lambda^+ y^2+\lambda x)} f(t)$, we would have $\sup_{s\in[0,\tau]} \hat{Z}_\lambda^+(s) \leq g_{\varepsilon_0}(\tau)$. Hence, defining $\hat{\mathcal{Z}}_\lambda^+(s) := \hat{Z}_\lambda^+(s) - \hat{f}(s)$, we have

$$\tilde{\mathbb{Q}}_\lambda^{x,y}\left(\mathbb{1}_{A_t^{\varepsilon,\delta}(\xi)} \tilde{\mathbb{Q}}_\lambda^{x,y}\left(\sup_{s\in[0,\tau]} \hat{Z}_\lambda^+(s) \leq g_{\varepsilon_0}(\tau)\Big|\tilde{\mathcal{G}}_\infty\right)\right)$$

$$\geq \tilde{\mathbb{Q}}_\lambda^{x,y}\left(\mathbb{1}_{A_t^{\varepsilon,\delta}(\xi)} \tilde{\mathbb{Q}}_\lambda^{x,y}\left(\sup_{s\in[0,\tau]} \hat{f}(s) \leq \frac{g_{\tilde{\varepsilon}_0}(\tau)}{2}; \sup_{s\in[0,\tau]} \hat{\mathcal{Z}}_\lambda^+(s) \leq \frac{g_{\varepsilon_0}(\tau)}{2}\Big|\tilde{\mathcal{G}}_\infty\right)\right)$$

$$= \tilde{\mathbb{Q}}_\lambda^{x,y}\left(\mathbb{1}_{A_t^{\varepsilon,\delta}(\xi)} \mathbb{1}_{\{\sup_{s\in[0,\tau]} \hat{f}(s)\leq g_{\tilde{\varepsilon}_0}(\tau)/2\}} \tilde{\mathbb{Q}}_\lambda^{x,y}\left(\sup_{s\in[0,\tau]} \hat{\mathcal{Z}}_\lambda^+(s) \leq \frac{g_{\varepsilon_0}(\tau)}{2}\Big|\tilde{\mathcal{G}}_\infty\right)\right)$$

since, conditional on $\tilde{\mathcal{G}}_\infty$, the supremum of $\hat{f}$ on $[0,\tau]$ is known.

We see from (50) that, *conditional on* $\tilde{\mathcal{G}}_\infty$, $\hat{\mathcal{Z}}_\lambda^+(t)$ is a submartingale. This is because the $\tilde{\mathbb{Q}}_\lambda^{x,y}$-conditional expectation of each of the $Z_\lambda^{(u)}$ in the sum

$$(51) \qquad e^{-(\psi_\lambda^+ y^2+\lambda x)} \sum_{u<\xi_t} e^{-E_\lambda^+ S_u} Z_\lambda^{(u)}(t-S_u)$$

is constant, so the expectation of the sum cannot decrease, and in fact this expectation increases every time there is a birth on the spine. Then by Doob's submartingale inequality we have

$$\tilde{\mathbb{Q}}^{x,y}_\lambda\left(\sup_{s\in[0,\tau]} \hat{\mathcal{Z}}^+_\lambda(s) \leq \frac{g_{\varepsilon_0}(\tau)}{2}\Big|\tilde{\mathcal{G}}_\infty\right)$$

$$= 1 - \tilde{\mathbb{Q}}^{x,y}_\lambda\left(\sup_{s\in[0,\tau]} \hat{\mathcal{Z}}^+_\lambda(s) \geq \frac{g_{\varepsilon_0}(\tau)}{2}\Big|\tilde{\mathcal{G}}_\infty\right)$$

$$\geq 1 - \frac{2}{g_{\varepsilon_0}(\tau)}\tilde{\mathbb{Q}}^{x,y}_\lambda(\hat{\mathcal{Z}}^+_\lambda(\tau)|\tilde{\mathcal{G}}_\infty).$$

We must note here that the expectation on the above line is not a priori finite. However, the expectation of each term in the sum (51) is bounded by $\sup_{s\in[0,\tau]} \hat{f}(s)$, which we have control over via an indicator function and so we do not have to worry about this expectation blowing up.

So we need to show that for all sufficiently large $\tau$ and all $y \in [y_0, y_1]$,

$$\tilde{\mathbb{Q}}^{x,y}_\lambda\left(\mathbb{1}_{A^{\varepsilon,\delta}_t(\xi)\cap\{\sup_{s\in[0,\tau]}\hat{f}(s)\leq g_{\tilde\varepsilon_0}(\tau)/2\}}\left(1 - \frac{2}{g_{\varepsilon_0}(\tau)}\tilde{\mathbb{Q}}^{x,y}_\lambda(\hat{\mathcal{Z}}^+_\lambda(\tau)|\tilde{\mathcal{G}}_\infty)\right)\right) > \varepsilon',$$

and hence also

$$\tilde{\mathbb{Q}}^{x,y}_\lambda\left(A^{\varepsilon,\delta}_t(\xi); \sup_{s\in[0,\tau]} \hat{Z}^+_\lambda(s) \leq g_{\varepsilon_0}(\tau)\right) > \varepsilon'$$

as required. This will follow by combining both parts of the following result.

LEMMA 12. *Fix $y_1 > y_0 > 0$ and $\varepsilon_0 > \tilde\varepsilon_0 > 0$.*

(i) *For all sufficiently small $\varepsilon, \delta > 0$, there exists some $\varepsilon' > 0$ and $\tilde T > 0$ such that for all $y \in [y_0, y_1]$ and all $t > \tilde T$,*

$$\tilde{\mathbb{Q}}^{x,y}_\lambda\left(A^{\varepsilon,\delta}_t(\xi); \sup_{s\in[0,\tau]} \hat{f}(s) \leq \frac{g_{\tilde\varepsilon_0}(\tau)}{2}\right) > \varepsilon'.$$

(ii) *As $\tau \to \infty$,*

$$\tilde{\mathbb{Q}}^{x,y}_\lambda\left(\frac{2}{g_{\varepsilon_0}(\tau)}\tilde{\mathbb{Q}}^{x,y}_\lambda(\hat{\mathcal{Z}}^+_\lambda(\tau)|\tilde{\mathcal{G}}_\infty); \sup_{s\in[0,\tau]} \hat{f}(s) \leq \frac{g_{\tilde\varepsilon_0}(\tau)}{2}\right) \to 0$$

*uniformly over $y \in [y_0, y_1]$.*

Then we have shown that, for any $\varepsilon_0 > 0$, $y_1 > y_0 > 0$, and sufficiently small $\varepsilon, \delta > 0$, there exists a $T > 0$ such that, for all $y \in [y_0, y_1]$ and all $t > T$,

$$t^{-1} \log P^{x,y}(A^{\varepsilon,\delta}_t) \geq -(\Theta(\beta,\kappa) + \varepsilon_0).$$

Finally, we observe that the probability $P^{x,y}(A_t^{\varepsilon,\delta})$ is trivially monotone increasing in both $\varepsilon$ and $\delta$, and so it follows that if the result is true for all sufficiently small $\varepsilon$ and $\delta$, it is in fact true for all $\varepsilon, \delta > 0$. This completes the proof of Theorem 7. □

PROOF OF LEMMA 12(i). We will prove Lemma 12(i) in a sequence of other lemmas, using a convenient coupling for the spine's type process.

First recall that, under $\tilde{\mathbb{Q}}_\lambda^{x,y}$, $\eta_s$ solves the SDE

$$d\eta_s = \sqrt{\theta}\, dB_s + \mu_\lambda \eta_s\, ds,$$

where $B_s$ is a $\tilde{\mathbb{Q}}_\lambda$-Brownian motion. Noting that $d(e^{-\mu_\lambda s}\eta_s) = e^{-\mu_\lambda s}\sqrt{\theta}\, dB_s$, we can construct $e^{-\mu_\lambda s}\eta_s$ as a time-change of a Brownian motion with

$$e^{-\mu_\lambda s}\eta_s - \eta_0 = \sqrt{\theta}\int_0^s e^{-\mu_\lambda w}\, dB_w = \sqrt{\frac{\theta}{2\mu_\lambda}}\tilde{B}(1 - e^{-2\mu_\lambda s}),$$

where $\tilde{B}$ is also a $\tilde{\mathbb{Q}}_\lambda^{x,y}$-Brownian motion started at the origin.

In this way, we will construct $\eta^y$ under $\mathbb{P}$ from a Brownian motion $B^y$ started at $y\sqrt{2\mu_\lambda/\theta}$ where, for $s \in [0, \infty)$,

$$\eta^y(s) = \sqrt{\frac{\theta}{2\mu_\lambda}}e^{\mu_\lambda s}B^y(1 - e^{-2\mu_\lambda s}).$$

To construct simultaneously all type processes $\eta^y$ under the same measure $\mathbb{P}$, we first construct the process $B^{y_0}$ as an independent Brownian motion started at $y_0\sqrt{2\mu_\lambda/\theta}$. Second, we construct the process $B^{y_1}$ by running an independent Brownian motion started at $y_1\sqrt{2\mu_\lambda/\theta}$ until it first hits the path of $B^{y_0}$, at which point we couple the two processes together. Next, for any other $y \in (y_0, y_1)$, we run an independent Brownian motion $B^y$ until it first meets with either the process $B^{y_0}$ below or $B^{y_1}$ above, at which point we couple it to the process it first hits.

Finally, we construct all the corresponding spatial processes $\xi^y$ under $\mathbb{P}$ from a single Brownian motion $W$ by defining

$$(52) \qquad \xi^y(s) = W\left(a\int_0^s \eta^y(w)^2\, dw\right) + \lambda a\int_0^s \eta^y(w)^2\, dw,$$

where $W$ is started at $x$ and is independent of the $B^y$ processes.

Constructed in this way, for each $y \in [y_0, y_1]$, the $\mathbb{P}$-law of $(\xi^y, \eta^y)$ is the same as the $\tilde{\mathbb{Q}}_\lambda^{x,y}$-law of $(\xi, \eta)$.

Fixing $\mu \in (0,1)$ and $K > \max\{y_1, 1\}$, we define the events and stopping times

$$A_\varepsilon^y := \left\{B^y(s) \in \left[1 - \frac{\varepsilon}{2\kappa}, 1 + \frac{\varepsilon}{2\kappa}\right], \forall s \in (1-\mu, 1]\right\},$$

$$T_0 := \inf\{t : B^{y_0}(t) = 0\}, \qquad T_K := \inf\{t : B^{y_1}(t) = K\},$$

$$\tilde{A}_{\varepsilon, K} := A_\varepsilon^{y_0} \cap A_\varepsilon^{y_1} \cap \{T_0 > 1\} \cap \{T_K > 1\}.$$

Then, clearly $\mathbb{P}(\tilde{A}_{\varepsilon,K}) > 0$ and, on the event $\tilde{A}_{\varepsilon,K}$, the coupling gives

$$0 < \eta^{y_0}(s) \leq \eta^{y}(s) \leq \eta^{y_1}(s) \leq K\sqrt{\frac{\theta}{2\mu_\lambda}}e^{\mu_\lambda s},$$

for all $s \geq 0$ and $y \in [y_0, y_1]$. Note that our construction also ensures that if event $A_\varepsilon^{y_0} \cap A_\varepsilon^{y_1}$ occurs then so must $A_\varepsilon^y$ for any $y \in [y_0, y_1]$, hence $A_\varepsilon^y \supset \tilde{A}_{\varepsilon,K}$.
□

LEMMA 13. *Let $\varepsilon > 0$. On event $\tilde{A}_{\varepsilon,K}$, there exists a deterministic time $s_0 = s_0(\varepsilon) > 0$ such that for all $\tau > s_0$,*

$$\sup_{s \in [0,\tau]} |\eta^y(s) - \bar{y}(s)| \leq \varepsilon\sqrt{t},$$

*for all $y \in [y_0, y_1]$.*

PROOF. Set $s_1 = -\frac{1}{2\mu_\lambda}\log\mu$ and then, on event $\tilde{A}_{\varepsilon,K}$, for all $\tau \geq s > s_1$ we have

$$\left|\eta^y(s) - \sqrt{\frac{\theta}{2\mu_\lambda}}e^{\mu_\lambda s}\right| \leq \frac{\varepsilon}{2\kappa}\sqrt{\frac{\theta}{2\mu_\lambda}}e^{\mu_\lambda s},$$

for all $y \in [y_0, y_1]$. Writing

$$\bar{y}(s) = \sqrt{\frac{\theta}{2\mu_\lambda}}e^{\mu_\lambda s}\left(\frac{1 - e^{-2\mu_\lambda s}}{1 - e^{-2\mu_\lambda \tau}}\right) \leq \sqrt{\frac{\theta}{2\mu_\lambda}}e^{\mu_\lambda s},$$

we see that there exists $s_2 = s_2(\varepsilon) > 0$ such that, for $\tau \geq s > s_2$,

$$\left|\bar{y}(s) - \sqrt{\frac{\theta}{2\mu_\lambda}}e^{\mu_\lambda s}\right| \leq \frac{\varepsilon}{2\kappa}\sqrt{\frac{\theta}{2\mu_\lambda}}e^{\mu_\lambda s}.$$

Taking $s_3(\varepsilon) = \max\{s_1, s_2(\varepsilon)\}$ now yields

(53) $$|\eta^y(s) - \bar{y}(s)| < \frac{\varepsilon}{\kappa}\sqrt{\frac{\theta}{2\mu_\lambda}}e^{\mu_\lambda s} \leq \varepsilon\sqrt{t}$$

for all $\tau \geq s > s_3$ and all $y \in [y_0, y_1]$.

Now consider $s \in [0, s_3]$. On $\tilde{A}_{\varepsilon,K}$ we have

$$|\eta^y(s) - \bar{y}(s)| \leq \sqrt{\frac{\theta}{2\mu_\lambda}}e^{\mu_\lambda s_3}(1 + K),$$

and hence for some $s_4(\varepsilon) > 0$ we have $|\eta^y(s) - \bar{y}(s)| \leq \varepsilon\sqrt{t}$ for all $\tau > s_4$, all $s \in [0, s_3]$, and all $y \in [y_0, y_1]$. Taking $s_0(\varepsilon) = \max\{s_3, s_4\}$ yields the result.

LEMMA 14. *Let $\delta > 0$. Then for all sufficiently small $\varepsilon$, there exists a deterministic $\tau_0 = \tau_0(\varepsilon, \delta) > 0$ such that, on $\tilde{A}_{\varepsilon,K}$, we have*

$$\text{(54)} \qquad \sup_{s \in [0,\tau]} \left| \int_0^s \eta^y(w)^2 \, dw - \int_0^s \bar{y}(w)^2 \, dw \right| < \delta t$$

*for all $\tau > \tau_0$ and all $y \in [y_0, y_1]$.*

PROOF. Given any $\delta > 0$, we first fix an $\varepsilon > 0$ sufficiently small such that $\varepsilon(2 + \frac{\varepsilon}{\kappa})\frac{\kappa}{2\mu_\lambda} < \frac{\delta}{4}$; this yields a corresponding $s_3 = s_3(\varepsilon)$, which is chosen as at equation (53). Given this $s_3$, we find $\tau_1 = \tau_1(\varepsilon, \delta) > 0$ such that, for all $\tau > \tau_1$,

$$(K^2 + 1) \int_0^{s_3} \frac{\theta}{2\mu_\lambda} e^{2\mu_\lambda w} \, dw < \frac{\delta}{4} t.$$

We now set $\tau_0 = \tau_0(\varepsilon, \delta) = \max\{s_3, \tau_1\}$. With this choice of $\varepsilon$ and $\tau_0$, we proceed to show that the inequality (54) is satisfied. Note that $\tau_0$ is deterministic and independent of $y$.

From equation (53) we see that, on $\tilde{A}_{\varepsilon,K}$ and for $s > s_3$,

$$\int_0^s \eta^y(w)^2 \, dw \geq \int_0^{s_3} \eta^y(w)^2 \, dw + \int_{s_3}^s \left( \bar{y}(w) - \frac{\varepsilon}{\kappa} \sqrt{\frac{\theta}{2\mu_\lambda}} e^{\mu_\lambda w} \right)^2 dw$$

$$\geq \int_0^s \bar{y}(w)^2 \, dw - \int_0^{s_3} \bar{y}(w)^2 \, dw - 2 \int_{s_3}^s \frac{\varepsilon\theta}{2\kappa\mu_\lambda} e^{2\mu_\lambda w} \, dw$$

$$\geq \int_0^s \bar{y}(w)^2 \, dw - \int_0^{s_3} \frac{\theta}{2\mu_\lambda} e^{2\mu_\lambda w} \, dw - (2\varepsilon)\frac{\kappa t}{2\mu_\lambda}$$

$$> \int_0^s \bar{y}(w)^2 \, dw - \frac{\delta}{2} t$$

for all $\tau > \tau_0$ and all $y \in [y_0, y_1]$. Similarly

$$\int_0^s \eta^y(w)^2 \, dw \leq \int_0^{s_3} \eta^y(w)^2 \, dw + \int_{s_3}^s \left( \bar{y}(w) + \frac{\varepsilon}{\kappa} \sqrt{\frac{\theta}{2\mu_\lambda}} e^{\mu_\lambda w} \right)^2 dw$$

$$\leq \int_0^s \bar{y}(w)^2 \, dw + \int_0^{s_3} \left( K \sqrt{\frac{\theta}{2\mu_\lambda}} e^{\mu_\lambda w} \right)^2 dw$$

$$+ \int_{s_3}^s \frac{\varepsilon}{\kappa} \left( 2 + \frac{\varepsilon}{\kappa} \right) \left( \sqrt{\frac{\theta}{2\mu_\lambda}} e^{\mu_\lambda w} \right)^2 dw$$

$$\leq \int_0^s \bar{y}(w)^2 \, dw + K^2 \int_0^{s_3} \frac{\theta}{2\mu_\lambda} e^{2\mu_\lambda w} \, dw + \varepsilon \left( 2 + \frac{\varepsilon}{\kappa} \right) \frac{\kappa t}{2\mu_\lambda}$$

$$< \int_0^s \bar{y}(w)^2 \, dw + \frac{\delta}{2} t$$

for all $\tau > \tau_0$ and all $y \in [y_0, y_1]$. Finally, for $s \in [0, s_3]$, on $\tilde{A}_{\varepsilon,K}$ we have

$$\left| \int_0^s \eta^y(w)^2 \, dw - \int_0^s \bar{y}(w)^2 \, dw \right| \leq \int_0^s \eta^y(w)^2 \, dw + \int_0^s \bar{y}(w)^2 \, dw$$

$$\leq (K^2 + 1) \int_0^{s_3} \frac{\theta}{2\mu_\lambda} e^{2\mu_\lambda w} \, dw < \delta t$$

for all $\tau > \tau_0$ and all $y \in [y_0, y_1]$. $\square$

LEMMA 15. *Let $\delta > 0$. Then for all sufficiently small $\varepsilon > 0$, there exists $\mathbb{P}$-almost everywhere on $\tilde{A}_{\varepsilon,K}$ a random time $S_0 = S_0(\delta, \varepsilon) < \infty$ such that*

$$\sup_{s \in [0,\tau]} \left| \xi^y(s) - \lambda a \int_0^s \eta^y(w)^2 \, dw \right| < \delta t,$$

*for all $y \in [y_0, y_1]$ and all $\tau > S_0$.*

PROOF. Given $\delta > 0$, choose any $\delta', \delta'' > 0$ such that $\delta'(|\beta/\lambda| + \delta'') < \delta$. Recalling the construction of $\xi^y$ at (52), we see from standard properties of Brownian motion that there almost surely exists some $S_1 = S_1(\delta') < \infty$ such that

$$\frac{1}{t} \sup_{s \in [0,t]} |W(s)| < \delta' \qquad \text{for all } t > S_1.$$

Then

(55) $$\sup_{s \in [0,\tau]} \left| W\left(a \int_0^s \eta^y(w)^2 \, dw\right) \right| < \delta' \left(a \int_0^\tau \eta^y(w)^2 \, dw\right)$$

for all $\tau$ such that $a \int_0^\tau \eta^y(w)^2 \, dw > S_1$, and by the coupling construction, on $\tilde{A}_{\varepsilon,K}$ this is true for all $y \in [y_0, y_1]$ if $a \int_0^\tau \eta^{y_0}(w)^2 \, dw > S_1$. Then there exists ($\mathbb{P}$-almost everywhere on $\tilde{A}_{\varepsilon,K}$) a random time $S_2 = S_2(\delta') < \infty$, which depends on $B^{y_0}$ and $S_1$, such that $a \int_0^\tau \eta^y(w)^2 \, dw > S_1$ for all $y \in [y_0, y_1]$ when $\tau > S_2$.

Now by Lemma 14, given $\delta''$ and a sufficiently small $\varepsilon$, there exists a deterministic $\tau_0 = \tau_0(\varepsilon, \delta'') > 0$ such that, on $\tilde{A}_{\varepsilon,K}$,

(56) $$a \int_0^\tau \eta^y(w)^2 \, dw \leq a \int_0^\tau \bar{y}(s)^2 \, ds + \delta'' t = \left(\left|\frac{\beta}{\lambda}\right| + \delta''\right) t$$

for all $\tau > \tau_0$ and all $y \in [y_0, y_1]$. Combining the inequalities at (55) and (56), we now see that, for $\tau > S_0 = S_0(\varepsilon, \delta', \delta'') = \max\{S_2, \tau_0\}$,

$$\sup_{s \in [0,\tau]} \left| \xi^y(s) - \lambda a \int_0^s \eta^y(w)^2 \, dw \right| = \sup_{s \in [0,\tau]} \left| W\left(a \int_0^s \eta^y(w)^2 \, dw\right) \right| < \delta t$$

for all $y \in [y_0, y_1]$. $\square$

On combining Lemmas 14 and 15 and recalling the definition of optimal path $\bar{x}$ at (37), we obtain the following:

LEMMA 16. *Let $\delta > 0$. Then for all sufficiently small $\varepsilon > 0$, there exists $\mathbb{P}$-almost everywhere on $\tilde{A}_{\varepsilon,K}$ a random time $\tilde{S}_0 = \tilde{S}_0(\delta, \varepsilon) < \infty$ such that*

$$\sup_{s \in [0,\tau]} |\xi^y(s) - \bar{x}(s)| < \delta t,$$

*for all $y \in [y_0, y_1]$, and all $\tau > \tilde{S}_0$.*

We may now draw everything together to finish the proof of Lemma 12(i). First we observe that since $\lambda < 0$, on event $A_t^{\varepsilon,\delta}(\xi)$,

$$\sup_{s \in [0,\tau]} e^{-(\psi_\lambda^+ y^2 + \lambda x)} \exp(\psi_\lambda^+ \eta_s^2 + \lambda \xi_s - E_\lambda^+ s)$$

$$\leq e^{-(\psi_\lambda^+ y^2 + \lambda x)} \exp(\psi_\lambda^+ (\kappa + \varepsilon)^2 t + \lambda(-\beta - \delta)t),$$

and so, given $\tilde{\varepsilon}_0$, we can choose first $\delta$ and then $\varepsilon$ sufficiently small so that

$$A_t^{\varepsilon,\delta}(\xi) \subset \left\{ \sup_{s \in [0,\tau]} \hat{f}(s) \leq \frac{g_{\tilde{\varepsilon}_0}(\tau)}{2} \right\}$$

and, from Lemmas 13 and 16, there exists a random time $\tilde{T} = \tilde{T}(\delta, \varepsilon) < \infty$ such that on $\tilde{A}_{\varepsilon,K}$ we have

$$\sup_{s \in [0,\tau]} |\eta^y(s) - \bar{y}(s)| < \varepsilon\sqrt{t}$$

and

$$\sup_{s \in [0,\tau]} |\xi^y(s) - \bar{x}(s)| < \delta t$$

for all $\tau > \tilde{T}$ and all $y \in [y_0, y_1]$. That is, $\tilde{A}_{\varepsilon,K} \cap \{\tilde{T} < \tau\} \subset A_t^{\varepsilon,\delta}(\xi^y)$ for each $y \in [y_0, y_1]$, with the slight abuse of notation that

$$A_t^{\varepsilon,\delta}(\xi^y) = \left\{ \sup_{s \in [0,\tau(t)]} |\eta^y(s) - \bar{y}(s)| < \varepsilon\sqrt{t};\ \sup_{s \in [0,\tau(t)]} |\xi^y(s) - \bar{x}(s)| < \delta t \right\}.$$

Note also that $\mathbb{P}(\tilde{A}_{\varepsilon,K}) > \varepsilon'$ for some $\varepsilon' > 0$.

Combining the above, for any $y \in [y_0, y_1]$ we have

$$\tilde{\mathbb{Q}}_\lambda^{x,y}\left(A_t^{\varepsilon,\delta}(\xi);\ \sup_{s \in [0,\tau]} \hat{f}(s) \leq \frac{g_{\tilde{\varepsilon}_0}(\tau)}{2}\right) = \tilde{\mathbb{Q}}_\lambda^{x,y}(A_t^{\varepsilon,\delta}(\xi)) = \mathbb{P}(A_t^{\varepsilon,\delta}(\xi^y))$$

$$\geq \mathbb{P}(\tilde{A}_{\varepsilon,K};\tilde{T} < \tau) \to \mathbb{P}(\tilde{A}_{\varepsilon,K})$$

as $\tau \to \infty$, as required. □

PROOF OF LEMMA 12(ii). Consider the expectation of the "sum term." We have

$$\tilde{\mathbb{Q}}_\lambda^{x,y}(\hat{\mathcal{Z}}_\lambda^+(\tau)|\tilde{\mathcal{G}}_\infty) = e^{-(\psi_\lambda^+ y^2 + \lambda x)}\tilde{\mathbb{Q}}_\lambda^{x,y}\left(\sum_{u<\xi_\tau} e^{-E_\lambda^+ S_u} Z_\lambda^{(u)}(t-S_u)\Big|\tilde{\mathcal{G}}_\infty\right)$$

$$= e^{-(\psi_\lambda^+ y^2 + \lambda x)} \sum_{u<\xi_\tau} e^{-E_\lambda^+ S_u} \tilde{\mathbb{Q}}_\lambda^{x,y}(Z_\lambda^{(u)}(t-S_u)|\tilde{\mathcal{G}}_\infty)$$

$$\leq e^{-(\psi_\lambda^+ y^2 + \lambda x)} n_\tau \max\{e^{\psi_\lambda^+ \eta(S_u)^2 + \lambda \xi(S_u) - E_\lambda^+ S_u} : u < \xi_\tau\}$$

$$\leq n_\tau \sup_{s\in[0,\tau]} \hat{f}(s).$$

Hence

(57) 
$$\tilde{\mathbb{Q}}_\lambda^{x,y}\left(\frac{2}{g_{\varepsilon_0}(\tau)}\tilde{\mathbb{Q}}_\lambda^{x,y}(\hat{\mathcal{Z}}_\lambda^+(\tau)|\tilde{\mathcal{G}}_\infty); \sup_{s\in[0,\tau]}\hat{f}(s) \leq \frac{g_{\tilde{\varepsilon}_0}(\tau)}{2}\right)$$

$$\leq \tilde{\mathbb{Q}}_\lambda^{x,y}\left(n_\tau \frac{g_{\tilde{\varepsilon}_0}(\tau)}{g_{\varepsilon_0}(\tau)}; \sup_{s\in[0,\tau]}\hat{f}(s) \leq \frac{g_{\tilde{\varepsilon}_0}(\tau)}{2}\right)$$

$$\leq e^{-(\varepsilon_0 - \tilde{\varepsilon}_0)t}\tilde{\mathbb{Q}}_\lambda^{x,y}(n_\tau),$$

and we can now calculate $\tilde{\mathbb{Q}}_\lambda^{x,y}(n_\tau) = \tilde{\mathbb{Q}}_\lambda^{x,y}(\tilde{\mathbb{Q}}_\lambda^{x,y}(n_\tau|\mathcal{G}_\infty))$, where $\mathcal{G}_\infty$ the $\sigma$-algebra generated by the path of the spine (*not* including the birth times). Conditional on $\mathcal{G}_\infty$, $n_\tau$ is a Poisson random variable with mean given by $\int_0^\tau 2(r\eta_s^2 + \rho)\,ds$, and using Fubini's theorem we have

$$\tilde{\mathbb{Q}}_\lambda^{x,y}\left(\int_0^\tau 2(r\eta_s^2 + \rho)\,ds\right)$$

$$= \int_0^\tau 2r\tilde{\mathbb{Q}}_\lambda^{x,y}(\eta_s^2)\,ds + 2\rho\tau$$

$$= \frac{r}{\mu_\lambda}\left(\frac{\theta}{2\mu_\lambda} + y^2\right)e^{2\mu_\lambda\tau} - \left(\frac{\theta}{2\mu_\lambda} + y^2\right)\frac{r}{\mu_\lambda\tau} - \frac{r\theta\tau}{\mu_\lambda} + 2\rho\tau$$

$$= \frac{r}{\mu_\lambda}\kappa^2 t + \frac{2y^2\kappa\mu_\lambda}{\theta}t + o(\tau).$$

So the $\tilde{\mathbb{Q}}_\lambda^{x,y}$-expectation of $n_\tau$ only grows linearly in $t$. Then since $\varepsilon_0 - \tilde{\varepsilon}_0 > 0$, the expression at (57) tends to 0 as $t \to \infty$. Moreover, the expectation at (57) is bounded by the $\tilde{\mathbb{Q}}_\lambda^{x,y_1}$-expectation, and hence the convergence is uniform over $y \in [y_0, y_1]$, as claimed. $\square$

**8. Martingale results.** In this section we recall some existing and prove some new martingale results that are intermediate steps in the proofs of Theorem 1 and the upper bound of Theorem 3. We recall from [13] that $E_\lambda^-$ [also written $E^-(\lambda)$] and $\Delta(\gamma)$ are Legendre conjugates with

$$(58) \qquad \Delta(\gamma) = \inf_{\lambda<0}\{E^-(\lambda) + \lambda\gamma\}, \qquad E^-(\lambda) = \sup_{\gamma>0}\{\Delta(\gamma) - \gamma\lambda\}.$$

If, for $\lambda_{\min} < \lambda < 0$, we write $\gamma_\lambda$ for the $\gamma$ value which achieves the supremum on the right-hand side of equation (58), then the functions $\lambda \mapsto \gamma_\lambda$ from $(-\lambda_{\min}, 0)$ to $(0, \infty)$, and $\gamma \mapsto \lambda_\gamma$ from $(0, \infty)$ to $(-\lambda_{\min}, 0)$ are inverses of each other and, of course, $\lambda_\gamma$ is the $\lambda$ value which achieves the infimum on the left-hand side of equation (58). In addition, we note that

$$(59) \qquad \gamma_\lambda = -\frac{\partial}{\partial \lambda} E^-(\lambda) = \sqrt{\frac{\theta a^2 \lambda^2}{\theta - 8r - 4a\lambda^2}},$$

that $E^-(\lambda)$ and $\Delta(\gamma)$ are convex functions, and that

$$(60) \qquad \begin{aligned} \tilde{c}(\theta) &= \sup\{\gamma : \Delta(\gamma) > 0\} = \inf\{-E^-(\lambda)/\lambda : \lambda_{\min} < \lambda < 0\} \\ &= \inf\{c_\lambda^- : \lambda_{\min} < \lambda < 0\} = c_{\tilde{\lambda}(\theta)}^-, \end{aligned}$$

where

$$(61) \qquad \begin{aligned} c_\lambda^- &:= -E_\lambda^-/\lambda \quad \text{and} \\ \tilde{\lambda}(\theta) &:= -\sqrt{\frac{2(\theta-8r)(\theta\rho + 2\rho^2 + r\theta)}{a(\theta+4\rho)^2}} \in (\lambda_{\min}, 0). \end{aligned}$$

A formula for $\tilde{c}(\theta)$ is given in equation (9). The following fundamental convergence result for the $Z_\lambda^-$ martingale was first partly proved in [13], but also see [9] for a more complete proof using "spine" techniques.

THEOREM 17. *Suppose* $\lambda \in (\lambda_{\min}, 0]$.

(i) *If* $\lambda \in (\tilde{\lambda}(\theta), 0]$, *the martingale* $Z_\lambda^-$ *is uniformly integrable and has an almost sure strictly positive limit.*

(ii) *If* $\lambda \leq \tilde{\lambda}(\theta)$, *then* $Z_\lambda^-(\infty) = 0$ *almost surely.*

The following convergence result was proved in [12] using martingales based on Hermite polynomials.

THEOREM 18. *Let* $\lambda \in (\tilde{\lambda}(\theta), 0]$ *and* $\alpha < 1/4$. *For each* $P^{x,y}$ *starting law and every continuous bounded function* $f : \mathbb{R} \mapsto \mathbb{R}$, *we have*

$$\sum_{u \in N_t} f(Y_u(t)) e^{\alpha Y_u(t)^2 + \lambda(X_u(t) + c_\lambda^- t)} \xrightarrow[\text{a.s.}]{} f_0 Z_\lambda^-(\infty),$$

*where*

$$(62) \qquad f_0 := \left(\frac{\mu_\lambda}{\theta}\right)^{1/4} \int_{\mathbb{R}} f(y) e^{\alpha y^2} e^{\psi_\lambda^- y^2} \phi(y)\, dy$$

*and $\phi(y)$ is the standard normal density.*

In this paper, we require a corollary to this theorem which specifies more precisely which particles contribute to the final limit.

COROLLARY 19. *Let $\lambda \in (\tilde{\lambda}(\theta), 0]$ and $\alpha < 1/4$. For each $P^{x,y}$ starting law and every continuous bounded function $f: \mathbb{R} \mapsto \mathbb{R}$, we have for every $\varepsilon > 0$*

$$(63) \qquad \sum_{u \in N_t} f(Y_u(t))\, e^{\alpha Y_u(t)^2 + \lambda X_u(t) - E_\lambda^- t} \mathbb{1}_{\{|X_u(t)/t + \gamma_\lambda| < \varepsilon\}} \xrightarrow[\text{a.s.}]{} f_0 Z_\lambda^-(\infty)$$

*where $\gamma_\lambda = -\frac{\partial}{\partial \lambda} E^-(\lambda)$ and $f_0$ is given at equation (62).*

This last result will enable us to show in Section 10 that the almost sure growth rate is at least as large as the expected growth rate, $D(\gamma) \geq \Delta(\gamma)$. It is easy to see from Corollary 19 that when $Z_\lambda^-(\infty) > 0$, there must exist at least one particle near to $-\gamma_\lambda t$ in space. Further, because of the decay rate of each term in the sum over particles at equation (63), it will be relatively straightforward to improve this to get the required exponential numbers of particles, $\exp(\Delta(\gamma)t)$, near $-\gamma_\lambda t$ for large times [as long as $Z_\lambda^-(\infty) > 0$].

The following result concerns the rate at which the martingales $Z_\lambda^+$ and $Z_\lambda^-$ converge to zero.

THEOREM 20. *Let $\lambda \in (\lambda_{\min}, 0)$. For every starting law, $P^{x,y}$,*

$$\frac{\log Z_\lambda^\pm(t)}{t} \to \lambda(c_\lambda^\pm - c_\lambda^*) \qquad a.s.$$

*where $c_\lambda^\pm$ is given at (5), and*

$$c_\lambda^* := \begin{cases} \tilde{c}(\theta), & \text{if } \lambda_{\min} < \lambda \leq \tilde{\lambda}(\theta), \\ c_\lambda^-, & \text{if } \tilde{\lambda}(\theta) \leq \lambda < 0. \end{cases}$$

COROLLARY 21. *If $\lambda \in (\lambda_{\min}, 0)$, then $Z_\lambda^+(t) \to 0$ $P^{x,y}$-almost surely.*

The rate of convergence of the $Z_\lambda^+$ martingale in part (i) of Theorem 20 is crucial in Section 9 to obtain the upper bound on the almost sure growth rate, $D(\gamma, \kappa) \leq \Delta(\gamma, \kappa)$. We also comment that *if* Corollary 19 were true for all $\alpha < \psi_\lambda^+$, then we could have gained this upper bound at that point.

Although Corollary 19 is only proven for $\alpha < 1/4$ (where we can utilize suitable Hermite expansions), we conjecture that it holds for all $\alpha < \psi_\lambda^+$.

PROOF OF COROLLARY 19. Let $\varepsilon > 0$ be small, $\mu := \lambda - \varepsilon$, $\lambda, \mu \in (\tilde{\lambda}(\theta), 0)$, $f$ be a positive, continuous bounded function, $\alpha < 1/4$ and note that $\gamma_\mu > \gamma_\lambda$. Then we have

$$\sum_{u \in N_t} f(Y_u(t)) e^{\alpha Y_u(t)^2 + \lambda X_u(t) - E_\lambda^- t} \mathbb{1}_{\{X_u(t) < -\gamma_\mu t\}}$$

$$\leq e^{(E_\mu^- - E_\lambda^- - \varepsilon \gamma_\mu)t} \sum_{u \in N_t} f(Y_u(t)) e^{\alpha Y_u(t)^2 + \mu X_u(t) - E_\mu^- t} \mathbb{1}_{\{X_u(t) < -\gamma_\mu t\}}$$

$$\leq \left( \sum_{u \in N_t} f(Y_u(t)) e^{\alpha Y_u(t)^2 + \mu X_u(t) - E_\mu^- t} \right) e^{-(E_\lambda^- - E_\mu^- + (\lambda - \mu)\gamma_\mu)t}.$$

Recall that $E^-(\lambda)$ is convex with $\frac{\partial^2}{\partial \lambda^2} E^-(\lambda) \geq 0$ and $\frac{\partial}{\partial \lambda} E^-(\lambda) = \gamma_\lambda$, so, from the Taylor expansion,

$$E_\lambda^- - E_\mu^- + (\mu - \lambda) \frac{\partial}{\partial \lambda} E^-(\lambda)$$
$$= \frac{(\mu - \lambda)^2}{2} \frac{\partial^2}{\partial \lambda^2} E^-(\lambda) + o((\mu - \lambda)^2).$$

Then taking $\varepsilon > 0$ small enough so that $E_\lambda^- - E_\mu^- + (\lambda - \mu)\gamma_\mu > 0$, and using Theorem 18, we find that for any $\delta > 0$

$$\limsup_{t \to \infty} \sum_{u \in N_t} f(Y_u(t)) e^{\alpha Y_u(t)^2 + \lambda X_u(t) - E_\lambda^- t} \mathbb{1}_{\{X_u(t) < -(\gamma_\lambda + \delta)t\}} = 0.$$

Similarly, we can show

$$\limsup_{t \to \infty} \sum_{u \in N_t} f(Y_u(t)) e^{\alpha Y_u(t)^2 + \lambda X_u(t) - E_\lambda^- t} \mathbb{1}_{\{X_u(t) > -(\gamma_\lambda - \delta)t\}} = 0,$$

and hence the only contribution to the limit comes from the particles near $-\gamma_\lambda t$ in space. Combining this with Theorem 18 we have

$$f_0 Z_\lambda^-(\infty) = \lim_{t \to \infty} \sum_{u \in N_t} f(Y_u(t)) e^{\alpha Y_u(t)^2 + \lambda X_u(t) - E_\lambda^- t}$$

$$= \lim_{t \to \infty} \sum_{u \in N_t} f(Y_u(t)) e^{\alpha Y_u(t)^2 + \lambda X_u(t) - E_\lambda^- t} \mathbb{1}_{\{|X_u(t)/t + \gamma_\lambda| < \delta\}}. \quad \square$$

PROOF OF THEOREM 20. We use a useful technique brought to our attention in [22]. Let $p \in (0, 1)$ so that, by Jensen's inequality, $Z_\lambda^\pm(t)^p$ is a supermartingale; then for $u, v > 0$ we have

$$(u + v)^p \leq u^p + v^p,$$

and hence

$$Z_\lambda^\pm(t)^p = \left| \sum_{u \in N_t} e^{\psi_\lambda^\pm Y_u(t)^2 + \lambda(X_u(t) + c_\lambda^\pm t)} \right|^p$$

$$\leq \left( \sum_{u \in N_t} e^{p\psi_\lambda^\pm Y_u(t)^2 + p\lambda(X_u(t) + c_{p\lambda}^- t)} \right) e^{p\lambda(c_\lambda^\pm - c_{p\lambda}^-)t}.$$

For any $\varepsilon > 0$, Doob's supermartingale inequality says

$$P\left( \sup_{s \leq w \leq s+t} Z_\lambda^\pm(w)^p > \varepsilon^p \right) \leq \frac{E Z_\lambda^\pm(s)^p}{\varepsilon^p}$$

$$\leq \varepsilon^{-p} \left( E \sum_{u \in N_s} e^{p\psi_\lambda^\pm Y_u(s)^2 + p\lambda(X_u(s) + c_{p\lambda}^- s)} \right) e^{p\lambda(c_\lambda^\pm - c_{p\lambda}^-)s},$$

and then

$$P\left( \sup_{s \leq w \leq s+t} e^{\delta w} Z_\lambda^\pm(w) > \varepsilon \right)$$

$$\leq P\left( \sup_{s \leq w \leq s+t} Z_\lambda^\pm(w)^p > e^{-p\delta(s+t)} \varepsilon^p \right)$$

$$\leq \varepsilon^{-p} e^{p\delta t} \left( E \sum_{u \in N_s} e^{p\psi_\lambda^\pm Y_u(s)^2 + p\lambda(X_u(s) + c_{p\lambda}^- s)} \right) e^{p(\lambda(c_\lambda^\pm - c_{p\lambda}^-) + \delta)s}.$$

Now, if we can choose $p \in (0,1)$ such that $\lambda(c_\lambda^\pm - c_{p\lambda}^-) + \delta < 0$ and $p\psi_\lambda^\pm < \psi_{p\lambda}^+$, we must have $e^{\delta u} Z_\lambda^\pm(u) \to 0$ almost surely by using a familiar Borel–Cantelli argument. [The condition $p\psi_\lambda^\pm < \psi_{p\lambda}^+$ guarantees that the expectation in the last line above tends to a finite limiting value, hence stays bounded over all times $s$, as can be checked by using formula (17), for example.]

For all $0 \leq p < 1$ we find $p\psi_\lambda^\pm < \psi_{p\lambda}^+$. Considering the graph of $c_\lambda^\pm$ we quickly see that, for $\lambda \in [\tilde{\lambda}(\theta), 0)$, taking $p$ as close to 1 as we like gives the best rate. For $\lambda \in [\lambda_{\min}, \tilde{\lambda}(\theta))$ we can choose $p$ so that $p\lambda = \tilde{\lambda}(\theta)$, which gives the best rate.

Recall from Theorem 17 that $Z_\lambda^-(\infty) > 0$ when $\lambda \in (\tilde{\lambda}(\theta), 0)$. Then, so far, we have proved the following:

LEMMA 22. *For every starting law, $P^{x,y}$, and for all $\varepsilon > 0$, if $\lambda \in (\lambda_{\min}, 0)$ then*

$$e^{-\varepsilon t} e^{-\lambda(c_\lambda^\pm - c_\lambda^*)t} Z_\lambda^\pm(t) \to 0 \qquad a.s.$$

*where*

$$c_\lambda^* := \begin{cases} \tilde{c}(\theta), & \text{if } \lambda_{\min} < \lambda \leq \tilde{\lambda}(\theta), \\ c_\lambda^-, & \text{if } \tilde{\lambda}(\theta) \leq \lambda < 0. \end{cases}$$

It is clear that this gives the required upper bound of

$$\limsup_{t \to \infty} \frac{\log Z_\lambda^\pm(t)}{t} \leq \lambda(c_\lambda^\pm - c_\lambda^*).$$

Now, for any $\varepsilon > 0$, if $\lambda \in (\lambda_{\min}, \tilde{\lambda}(\theta)]$ then

$$e^{\varepsilon t} \sum_{u \in N_t} e^{\psi_\lambda^\pm Y_u(t)^2 + \lambda(X_u(t) + \tilde{c}(\theta)t)} \geq e^{\lambda(L_t + \tilde{c}(\theta)t) + \varepsilon t} \to \infty \quad \text{a.s.}$$

since we know that $L_t := \inf\{X_u(t) : u \in N(t)\}$ satisfies $L_t/t \to -\tilde{c}(\theta)$ a.s. Otherwise, with $\lambda \in (\tilde{\lambda}(\theta), 0)$,

$$e^{\varepsilon t} \sum_{u \in N_t} e^{\psi_\lambda^+ Y_u(t)^2 + \lambda(X_u(t) + c_\lambda^- t)} \geq e^{\varepsilon t} Z_\lambda^-(t) \to \infty \quad \text{a.s.}$$

since here $Z_\lambda^-(\infty) > 0$ a.s. Thus, in all cases,

$$\liminf_{t \to \infty} \frac{\log Z_\lambda^\pm(t)}{t} \geq \lambda(c_\lambda^\pm - c_\lambda^*),$$

which completes the proof of Theorem 20. □

**9. Proof of Theorem 3. Upper bound.** The idea for the upper bound proof is to overestimate indicator function by exponentials, and then rearrange the expressions to form martingale terms.

Simply observe that for $\lambda \in (\lambda_{\min}, 0)$,

$$N_t(\gamma, [\kappa\sqrt{t}, \infty)) = \sum_{u \in N_t} \mathbb{1}_{\{X_u(t) \leq -\gamma t; Y_u(t) \geq \kappa\sqrt{t}\}}$$

$$\leq \sum_{u \in N_t} \mathbb{1}_{\{X_u(t) \leq -\gamma t; Y_u(t)^2 \geq \kappa^2 t\}} e^{\psi_\lambda^+(Y_u(t)^2 - \kappa^2 t) + \lambda(X_u(t) + \gamma t)}$$

(64)
$$\leq e^{(E_\lambda^+ - \kappa^2 \psi_\lambda^+ + \lambda\gamma)t} \sum_{u \in N_t} e^{\psi_\lambda^+ Y_u(t)^2 + \lambda X_u(t) - E_\lambda^+ t}$$

$$\leq e^{-\lambda(c_\lambda^+ - c_\lambda^-)t} Z_\lambda^+(t) e^{(E_\lambda^- + \lambda\gamma - \kappa^2 \psi_\lambda^+)t},$$

where $E_\lambda^\pm = -\lambda c_\lambda^\pm$.

Recall from equations (11) and (32) that $E_\lambda^- + \lambda\gamma - \kappa^2 \psi_\lambda^+$ has a minimal value of $\Delta(\gamma, \kappa)$ achieved when $\lambda = \bar{\lambda}(\gamma, \kappa)$. Since $\tilde{c}(\theta)$ is the minimal value of $c_\lambda$, Theorem 20 implies that

(65) $$\limsup_{t \to \infty} t^{-1} \log Z_\lambda^+(t) \leq \lambda(c_\lambda^+ - c_\lambda^-)$$

almost surely for all $\lambda \in (\lambda_{\min}, 0)$.

In cases where $\Delta(\gamma, \kappa) < 0$, we can use the optimal value for $\lambda$, Theorem 20 and trivially note that $N_t(\gamma, [\kappa\sqrt{t}, \infty))$ is integer valued to deduce that

$$\sum_{u \in N_t} \mathbb{1}_{\{Y_u(t) \geq \kappa\sqrt{t}; X_u(t) \leq -\gamma t\}} = 0$$

eventually, almost surely. Hence, $D(\gamma, \kappa) = -\infty$ almost surely if $\Delta(\gamma, \kappa) < 0$.

Otherwise we have $\Delta(\gamma, \kappa) \geq 0$, which in fact guarantees that $\gamma \in (0, \tilde{c}(\theta)]$ and hence $\bar{\lambda}(\gamma, \kappa) \in [\tilde{\lambda}(\theta), 0)$. Then since

$$\limsup_{t \to \infty} t^{-1} \log N_t(\gamma, [\kappa\sqrt{t}, \infty))$$
$$\leq \limsup_{t \to \infty} t^{-1} \log(e^{-\lambda(c_\lambda^+ - c_\lambda^-)t} Z_\lambda^+(t)) + (E_\lambda^- + \lambda\gamma - \kappa^2 \psi_\lambda^+)$$

we can again make use of Theorem 20 and the minimizing $\lambda$ value, $\bar{\lambda}(\gamma, \kappa)$, to get the bound

$$\limsup_{t \to \infty} t^{-1} \log N_t(\gamma, [\kappa\sqrt{t}, \infty)) \leq \Delta(\gamma, \kappa) \qquad \text{almost surely},$$

as desired.

Notice that, when $\Delta(\gamma, \kappa) = 0$, the right-hand side of the inequality at (64) will tend to infinity (see Corollary 19). Then, on the boundary, we have only shown that $\limsup t^{-1} \log N_t(\gamma, [\kappa\sqrt{t}, \infty)) \leq 0$.

**10. Proof of Theorem 1. The spatial growth rate.** We first bound the spatial growth rate above. Suppose that $C \subset \mathbb{R}$ is Borel-measurable with $\int_C e^{\psi_\lambda^- y^2} \phi(y) \, dy > 0$. Let $\lambda \in (\lambda_{\min}, 0)$, then

$$\sum_{u \in N_t} \mathbb{1}_{\{X_u(t) \leq -\gamma t; Y_u(t) \in C\}} \leq \sum_{u \in N_t} \mathbb{1}_{\{Y_u(t) \in C\}} e^{\lambda(X_u(t) + \gamma t)}$$
$$= e^{(E_\lambda^- + \lambda\gamma)t} \sum_{u \in N_t} \mathbb{1}_{\{Y_u(t) \in C\}} e^{\lambda X_u(t) - E_\lambda^- t}$$
$$\leq e^{(E_\lambda^- + \lambda\gamma)t} Z_\lambda^-(t).$$

Recalling equations (8) and (19), we therefore have

$$\sum_{u \in N_t} \mathbb{1}_{\{X_u(t) \leq -\gamma t; Y_u(t) \in C\}} \leq e^{\Delta(\gamma)t} Z_{\lambda_\gamma}^-(t).$$

Now if $\gamma \geq \tilde{c}(\theta)$, corresponding to $\lambda_\gamma \in (\lambda_{\min}, \tilde{\lambda}(\theta)]$ and having $\Delta(\gamma) \leq 0$, we know from Theorem 17 that $Z_{\lambda_\gamma}^-(\infty) = 0$ almost surely. Then,

$$\gamma > \tilde{c}(\theta) \quad \Rightarrow \quad \sum_{u \in N_t} \mathbb{1}_{\{X_u(t) \leq -\gamma t; Y_u(t) \in C\}} = 0 \qquad \text{eventually, a.s.}$$

Otherwise, if $\gamma \in (0, \tilde{c}(\theta))$, corresponding to $\lambda_\gamma \in (\tilde{\lambda}(\theta), 0)$ and having $\Delta(\gamma) > 0$, Theorem 17 tells us that $Z^-_{\lambda_\gamma}(\infty) > 0$ almost surely, hence

$$\limsup_{t \to \infty} t^{-1} \log \sum_{u \in N_t} \mathbb{1}_{\{X_u(t) \leq -\gamma t; Y_u(t) \in C\}} \leq \Delta(\gamma).$$

Now we bound the growth rate from below. Let $\varepsilon > 0$ be small, $\tilde{\lambda}(\theta) < \lambda < 0$, and $\mu = \lambda - \varepsilon$. We recall now that $E^-_\lambda$ is convex so $\frac{\partial^2 E^-_\lambda}{\partial \lambda^2} \geq 0$ and $\gamma_\mu > \gamma_\lambda$. Then

$$\sum_{u \in N_t} e^{\lambda X_u(t) - E^-_\lambda t} \mathbb{1}_{\{-(\gamma_\lambda + \varepsilon)t \leq X_u(t) \leq -(\gamma_\lambda - \varepsilon)t; Y_u(t) \in C\}}$$

$$\leq \sum_{u \in N_t} e^{\lambda(-(\gamma_\lambda + \varepsilon)t) - E^-_\lambda t} \mathbb{1}_{\{-(\gamma_\lambda + \varepsilon)t \leq X_u(t) \leq -(\gamma_\lambda - \varepsilon)t; Y_u(t) \in C\}}$$

$$= e^{(-\lambda \gamma_\lambda - E^-_\lambda - \lambda \varepsilon)t} \sum_{u \in N_t} \mathbb{1}_{\{-(\gamma_\lambda + \varepsilon)t \leq X_u(t) \leq -(\gamma_\lambda - \varepsilon)t; Y_u(t) \in C\}}$$

$$\leq e^{(-\lambda \gamma_\lambda - E^-_\lambda - \lambda \varepsilon)t} \sum_{u \in N_t} \mathbb{1}_{\{X_u(t) \leq -(\gamma_\lambda - \varepsilon)t; Y_u(t) \in C\}}.$$

Then

$$t^{-1} \log \sum_{u \in N_t} \mathbb{1}_{\{Y_u(t) \in C\}} e^{\lambda X_u(t) - E^-_\lambda t} \mathbb{1}_{\{|X_u(t)/t + \gamma_\lambda| < \varepsilon\}}$$

$$\leq -\lambda \gamma_\lambda - E^-_\lambda - \lambda \varepsilon + t^{-1} \log \sum_{u \in N_t} \mathbb{1}_{\{X_u(t) \leq -(\gamma_\lambda - \varepsilon)t; Y_u(t) \in C\}}.$$

Letting $t \to \infty$, using Corollary 19 and remembering that for $\tilde{\lambda}(\theta) < \lambda \leq 0$ we have $Z^-_\lambda(\infty) > 0$ a.s., we find

$$0 \leq -\lambda \gamma_\lambda - E^-_\lambda - \lambda \varepsilon + \liminf_{t \to \infty} t^{-1} \log \sum_{u \in N_t} \mathbb{1}_{\{X_u(t) \leq -(\gamma_\lambda - \varepsilon)t; Y_u(t) \in C\}}$$

and as $\varepsilon > 0$ can be arbitrarily small we have

$$\liminf_{t \to \infty} t^{-1} \log \sum_{u \in N_t} \mathbb{1}_{\{X_u(t) \leq -\gamma_\lambda t; Y_u(t) \in C\}} \geq E^-_\lambda + \lambda \gamma_\lambda.$$

Equivalently,

$$\liminf_{t \to \infty} t^{-1} \log \sum_{u \in N_t} \mathbb{1}_{\{X_u(t) \leq -\gamma t; Y_u(t) \in C\}} \geq E^-_{\lambda_\gamma} + \lambda_\gamma \gamma = \Delta(\gamma)$$

and hence the lim sup and lim inf agree as required.

We note that these proofs will easily adapt to cover a multi-type branching Brownian motion where the types evolve as a finite state Markov chain, such as found in [2], where it will also be possible to prove the analogous convergence theorem required when we have a finite type space by adapting the proof of Theorem 18 found in [12].

In the standard branching Brownian motion case things are even simpler to adapt (where, of course, there is no need for any convergence result akin to Theorem 18). All the information necessary is contained in the martingales $\sum_{u \in N_t} \exp(\lambda X_u(t) - (\lambda^2/2 + r)t)$ studied by Neveu [22] and, as first came to our attention during discussions with J. Warren, the martingale with parameter $\lambda$ can *only* be capable of "counting" particles near $\gamma_\lambda t$ in space at large times $t$, so when this martingale is uniformly integrable particles *must* perpetually be found with the corresponding speed. Of course, in this case more precise results, in the spirit of Watanabe [25], also exist.

**Acknowledgments.** We would like to thank two anonymous referees for providing extremely helpful and thorough reviews of earlier incarnations of this manuscript. Their numerous invaluable comments led to a much improved presentation of this work.

Y. GIT
STATISTICAL LABORATORY
CAMBRIDGE UNIVERSITY
22 MILL STREET
CAMBRIDGE CB1 2HP
UK
E-MAIL: Yoav.Git@gmail.com

J. W. HARRIS
DEPARTMENT OF MATHEMATICS
UNIVERSITY OF BRISTOL
UNIVERSITY WALK
BRISTOL BS8 1TW
UK
E-MAIL: john.harris@bristol.ac.uk

S. C. HARRIS
DEPARTMENT OF MATHEMATICAL SCIENCES
UNIVERSITY OF BATH
BATH BA2 7AY
UK
E-MAIL: s.c.harris@bath.ac.uk
URL: http://people.bath.ac.uk/massch/